\documentclass[a4paper, 12pt]{article}
\usepackage[a4paper,left=1.45in,right=1.45in,top=1.45in,bottom=1.45in]{geometry}
\usepackage{amsmath, amsfonts, amsthm, MnSymbol, slashed, booktabs}
\usepackage[bookmarks=true, breaklinks=false, colorlinks=true, citecolor=blue, bookmarksnumbered=true]{hyperref}
\usepackage[numbers]{natbib}  
\usepackage{relsize}
\usepackage{mathrsfs}
\usepackage{enumitem}
 \usepackage[]{units}
\usepackage{tikz}
\usetikzlibrary{arrows,positioning,matrix,arrows}

\renewenvironment{quote}
               {\list{}{\setlength\leftmargin{1.5\parindent}
               \setlength\rightmargin{\leftmargin}} %
                \item\relax}
               {\endlist\ignorespaces}

\newcounter{ncount}                 
\newcounter{ncounts}

\newcommand\functionto{\longrightarrow}
\newcommand\lonlyif{\rightarrow}
\newcommand\liff{\leftrightarrow}
\newcommand\isomorphic{\cong}

\makeatletter
\def\s@btitle{\relax}
\def\subtitle#1{\gdef\s@btitle{#1}}
\def\em@iladdress{\relax}
\def\emailaddress#1{\gdef\em@iladdress{#1}}
\def\inst@tution{\relax}
\def\institution#1{\gdef\inst@tution{#1}}

\newcommand\maketitletop{%
	\begingroup%
	\begin{center}
		{\Large\bfseries \@title \par}%
		\if\s@btitle\relax%
		\else\typeout{[subtitle]}{\large\bfseries \s@btitle}%
		\fi
		\vskip 1em%
		{\@author \par}%
		\if\em@iladdress\relax%
		\else\typeout{[emailaddress]}{\em@iladdress \par}%
		\fi
		\if\inst@tution\relax%
		\else\typeout{[instutiton]}{\inst@tution}%
		\fi
		\vskip 1em%
		\@date \par%
	\end{center}\par
	\endgroup}
\makeatother

\newtheoremstyle{generic}
	{1em} 
	{1em} 
	{} 
	{}
	{\normalfont }
	{. }
	{0em}
	{}
	
\newcounter{dummy}\numberwithin{dummy}{section}
\theoremstyle{generic}
\newtheorem{thm}[dummy]{\normalfont\bfseries{Theorem}}
\newtheorem{define}[dummy]{\normalfont\bfseries{Definition}}

\newtheorem{lem}[dummy]{\normalfont\bfseries{Lemma}}

\newcommand\model[1]{\mathcal{#1}}
\newcommand\lang[1]{\mathscr{#1}}

\newcommand{\domain}{\text{dom}}
\newcommand{\range}{\text{rng}}
\newcommand{\yesstroke}{^+}
\newcommand{\nostroke}{^-}
\newcommand{\langyes}{\lang{L}\yesstroke}
\newcommand{\langno}{\lang{L}\nostroke}
\newcommand{\langmaybe}{\lang{L}}
\newcommand{\precyes}{\prec\yesstroke}
\newcommand{\precno}{\prec\nostroke}
\newcommand{\succyes}{\succ\yesstroke}
\newcommand{\succno}{\succ\nostroke}

\newcommand{\indiscernibleyes}{=}
\newcommand{\indiscernibleno}{\approx} 
\newcommand{\nindiscernibleyes}{\neq}
\newcommand{\nindiscernibleno}{\napprox}
\newcommand{\symmetric}{\mathrel{\textsf{s}}}
\newcommand{\nsymmetric}{\mathrel{\slashed{\symmetric}}}
\newcommand{\relativeto}{\mathrel{\textsf{r}}}
\newcommand{\nrelativeto}{\mathrel{\slashed{\relativeto}}}

\newcommand{\relstroke}{_\textsf{p}}
\newcommand{\monadstroke}{_\textsf{m}}
\newcommand{\totalstroke}{}
\newcommand{\pairstroke}{_\textsf{p}}
\newcommand{\barestroke}{_\textsf{m}}

\newcommand{\yesmonad}{\mathrel{\indiscernibleyes\monadstroke}}
\newcommand{\nyesmonad}{\mathrel{\nindiscernibleyes\monadstroke}}

\newcommand{\yesrel}{\mathrel{\indiscernibleyes\relstroke}}

\newcommand{\foi}{\indiscernibleno}
\newcommand{\nfoi}{\nindiscernibleno}
\newcommand{\nomonad}{\mathrel{\indiscernibleno\monadstroke}}

\newcommand{\norel}{\mathrel{\indiscernibleno\relstroke}}
\newcommand{\nnorel}{\mathrel{\nindiscernibleno\relstroke}}

\newcommand{\symtotal}{\symmetric\totalstroke}
\newcommand{\sympair}{\symmetric\pairstroke}
\newcommand{\symbare}{\symmetric\barestroke}

\newcommand{\nsymbare}{\nsymmetric\barestroke}

\newcommand{\relativetotal}{\relativeto\totalstroke}
\newcommand{\relativepair}{\relativeto\pairstroke}
\newcommand{\relativebare}{\relativeto\barestroke}
\newcommand{\nrelativetotal}{\nrelativeto\totalstroke}

\newcommand{\nrelativebare}{\nrelativeto\barestroke}

\newcommand{\verified}{$\checkmark$}
\newcommand{\falsified}{$\times$}
\newcommand{\finitecase}{\emph{\textbf{f}}}
\newcommand{\finitequotient}{\emph{\textbf{fq}}}

\newcommand\quotient[1]{\underline{#1}}

\newcommand{\Blank}{A}
\newcommand{\Black}{B}
\newcommand{\ThreeSphere}{C}
\newcommand{\Rotational}{D}
\newcommand{\TwoBallBig}{E}
\newcommand{\NearCorrespondence}{F}
\newcommand{\Hexagon}{G}
\newcommand{\NotDefineIdentity}{H}
\newcommand{\ThreeThreeSphere}{I}
\newcommand{\ThreeBallBig}{J}
\newcommand{\TripleRotate}{K}

\newcommand{\GaloisNear}{\textbf{N}}
\newcommand{\GaloisIso}{\textbf{I}}
\newcommand{\galup}{\textbf{e}}
\newcommand{\galdown}{\textbf{c}}

\newcommand{\eldiag}{\text{Diag}}

\newcommand{\edyes}{\eldiag\yesstroke}

\title{Grades of Discrimination}
\subtitle{Indiscernibility, symmetry, and relativity}
\author{Tim Button}
\emailaddress{button@cantab.net}
\institution{University of Cambridge}
\date{forthcoming in \emph{Notre Dame Journal of Formal Logic}\\final version submitted 26.iii.2015}
\begin{document}
 \maketitletop

\noindent

\section{Introduction}\label{s:Preliminaries}
There are several relations which may fall short of genuine identity, but which behave like identity in important respects. Such \emph{grades of discrimination} have recently been the subject of much philosophical and technical discussion.

Much of this discussion has been fuelled by considering \emph{the Principle of the Identity of Indiscernibles}: the claim that indiscernible objects are always identical. The Principle is obviously of direct metaphysical interest  (see \cite{Hawley:II}). But, within the philosophy of mathematics, the Principle has risen to prominence via the question of whether platonistically-minded structuralists can countenance structures with indiscernible but distinct positions (see \cite[\S1]{Shapiro:STUR}). And, within the philosophy of physics, the central question has been whether quantum mechanics presented real-world counterexamples to the Principle (see \cite{Muller:RR}). As discussion has progressed, though, it has become increasingly clear that we must distinguish between different versions of `the' Principle, corresponding to different notions of indiscernibility. This has spurred several philosophers to investigate the logical properties of these different notions (see \cite{CaultonButterfield:KILM}, \cite{Ketland:II},  \cite{LadymanLinneboPettigrew:IDPL}). 

This paper completes that logical investigation. It exhaustively details, not just the properties of grades of indiscernibility, but the properties of all of the grades of {discrimination}. Indeed, this paper answers all of the mathematical questions that are natural at this level of abstraction. 

There are three broad families of grades of discrimination. Grades of \emph{indiscernibility} are defined in terms of satisfaction of certain first-order formulas, either with or without access to a primitive symbol that stands for genuine identity. They have been the focus of much recent philosophical attention. Grades of \emph{symmetry} are defined in terms of isomorphisms. More specifically, they are defined in terms of symmetries (also known as automorphisms) on a structure. These grades have received some philosophical attention, though in a slightly less cohesive way than the grades of indiscernibility.  Finally, grades of \emph{relativity} are defined in terms of relativeness correspondences, analogously to the grades of symmetry. The notion of a relativeness correspondence has been studied by model-theorists, but is entirely absent from the philosophical literature on grades of discrimination. This paper rectifies this situation, introducing grades of relativity for the first time. 

I mentioned, earlier, that the Principle of the Identity of Indiscernibles has been the main motivating force for interest in grades of discrimination. But it is now worth pausing to consider broader reasons for investigating the logical properties of the grades of discrimination. 

The simplest reason to care about grades of discrimination is that they allow us to calibrate relationships of similarity and difference. More ambitiously, though, we might hope that some grade of discrimination will provide us with a genuinely illuminating answer to the question: {When are objects identical?} To take a simple example: set theory tells us that sets are identical iff they share all their members. Consequently, some grade of indiscernibility provides a suitable criterion of identity in set-theoretic contexts. To take a more contentious example: we might somehow become convinced that nature abhors a (non-trivial) symmetry. If so, then some grade of symmetry will provide a suitable criterion of identity in empirical contexts. The general hope, then, is that our grades of discrimination may furnish us with some non-trivial \emph{criterion of identity} (in some context or other). 

This search for a non-trivial criterion of identity need not be \emph{reductive}. We might simply seek an illuminating {constraint} upon the conditions under which objects can be distinct. That said, some philosophers have hoped to find a {reductive} criterion for identity; that is, they have hoped to replace the identity primitive with some defined grade of discrimination. This reductive ambition is most prominent among those who have defended some Principle of Identity of Indiscernibles; such philosophers have therefore focussed on the various grades of indiscernibility. However, reductive ambitions might, in principle, be served equally well by considering either grades of symmetry or grades of relativity. (I revisit this in \S\ref{s:Capturing}.)

Advancing a criterion of identity is not, however, simply a matter of selecting some grade of discrimination. As we shall see, each grade of discrimination is defined with respect to a (model-theoretic) signature. So, consider a signature which contains just a few monadic predicates which stand for eye colour. If we present a non-trivial criterion of identity, in the form of a grade of discrimination defined with respect to \emph{this} signature, then we shall be forced to say, absurdly, that there is at most one person with brown eyes. Consequently, any philosopher who wants to advance a non-trivial criterion of identity must not only select some appropriate grade of discrimination, but must also stipulate the particular signature she has in mind. (I shall revisit this point several times below.)

In this paper, though, I am not aiming to \emph{advance} any particular criterion of identity. My aim is only to provide a mathematical toolkit for anyone who is interested in criteria of identity, whether reductive or non-reductive. 

That toolkit is structured around three main results. Theorem \ref{thm:BigMap} completely characterises the entailments between the grades of discrimination. Theorem \ref{thm:GaloisConnection} establishes a Galois Connection between isomorphisms and relativeness correspondences, which enables us better to understand the relationships between grades of symmetry and grades of relativity. And Theorem \ref{thm:NoSven}, is a Beth--Svenonius Theorem for logics without identity. By combining these three results, I answer several subsidiary questions concerning the grades of discrimination, including: which grades are equivalence relations (\S\ref{s:Equivalence}); which grades can be captured using sets of first-order formulas (\S\ref{s:Capturing}); how they behave in finitary cases (\S\ref{s:FinitaryCase}); and how they behave in elementary extensions of structures (\S\ref{s:DefinabilityTheory} and \S\ref{s:AllElementaryExtensions}).

I now state some notational conventions.

I always use `$\langmaybe$' to denote an arbitrary signature, i.e.\ a collection of constants, predicates and function-symbols. The philosophical discussion of grades of indiscernibility tends to be restricted to \emph{relational} signatures, i.e.\ signatures which contain only predicates.\footnote{An exception is \citeauthor{LadymanLinneboPettigrew:IDPL} \cite[\S6]{LadymanLinneboPettigrew:IDPL}.} There are reasonable philosophical motivations for this: if we assume that each constant names exactly \emph{one} object, then we seem to presuppose that we understand rather a lot about the notion of identity before we begin (see e.g.\ \cite{Black:II}, \cite[pp.\ 40--1]{CaultonButterfield:KILM}); more generally, the very idea of a \emph{function} seems to presuppose the notion of identity; hence, if we want to avoid prejudging certain philosophical questions about identity, it might be wise to restrict our attention to relational signatures. The model-theoretic discussion of these issues is, though, less often restricted to relational signatures. There is a sensible technical motivation for this: many of the results hold in the more general case. Since this paper aims to provide philosophers with technical results, I shall  allow signatures to contain both constants and functions, but I shall comment on the relational case when it is interestingly different. For technical ease, I treat constants as $0$-place function-symbols.

Where $\langmaybe$ is a signature, the $\langyes$-formulas are the first-order formulas  formed in the usual way using any $\langmaybe$-symbols and any symbols from standard first-order logic with identity. In particular, then, they may contain the symbol `$=$', which always stands for (genuine) identity. The $\langno$-formulas are those formed \emph{without} using the symbol `$=$'. $\langyes_n$ is the set of $\langyes$-formulas with free variables among `$v_1$', $\ldots$, `$v_n$'; similarly for $\langno_n$. 

I use swash fonts for structures and italic fonts for their associated domains. So, where $\model{M}$ is an $\langmaybe$-structure, its domain is $M$. Where $\overline{e} = \langle e_1, \ldots, e_n\rangle$ and $\pi$ is a function, $\pi(\overline{e}) = \langle \pi(e_1), \ldots, \pi(e_n)\rangle$. Where $\Pi$ is a two place relation, I write $\overline{d}\Pi\overline{e}$ to abbreviate $\langle d_{1}, e_{1}\rangle, \ldots, \langle d_{n}, e_{n}\rangle \in \Pi$.

\section{Twelve grades of discrimination}\label{s:TwelveGrades}
I start by defining six \emph{grades of indiscernibility}; three grades of $\langno$-indiscerni\-bi\-lity, and three grades of $\langyes$-indiscernibility.\footnote{This family of definitions has a long philosophical heritage, e.g.: \citeauthor{HilbertBernays:GM} \cite[\S5]{HilbertBernays:GM}; Quine \cite[pp.\ 230--2]{Quine:WO}, \cite{Quine:GD}; \citeauthor{CaultonButterfield:KILM} \cite[\S2.1, \S3.2]{CaultonButterfield:KILM}; \citeauthor{Ketland:II} \cite[pp.\ 306--7]{Ketland:SII}, \cite[Definitions 2.3, 2.5]{Ketland:II};  \citeauthor{LadymanLinneboPettigrew:IDPL} \cite[Definition 3.1, \S6.4]{LadymanLinneboPettigrew:IDPL}.}
	\begin{define}\label{def:Indiscernibility}
		For any $\langmaybe$-structure $\model{M}$ with $a, b \in M$:
	\begin{enumerate}[nolistsep]
		\item $a \nomonad b$ in $\model{M}$ \text{iff} $\model{M} \models \phi(a) \liff \phi(b)$ for all $\phi \in \langno_1$
		\item $a \norel b$ in $\model{M}$ \text{iff} $\model{M} \models \phi(a,b) \liff \phi(b,a)$ for all $\phi \in \langno_2$
		\item $a \foi b$ in $\model{M}$ \text{iff} $\model{M} \models \phi(a, \overline{e}) \liff \phi(b, \overline{e})$ for all $n < \omega$, all $\phi \in \langno_{n+1}$ and all $\overline{e} \in M^n$
	\end{enumerate}	
	Similarly:
	\begin{enumerate}[nolistsep]\addtocounter{enumi}{3}
		\item $a \yesmonad b$ in $\model{M}$ \text{iff} $\model{M} \models \phi(a) \liff \phi(b)$ for all $\phi \in \langyes_1$ 
		\item $a \yesrel b$ in $\model{M}$ \text{iff} $\model{M} \models \phi(a,b) \liff \phi(b,a)$ for all $\phi \in \langyes_2$
		\item $a = b$ in $\model{M}$ \text{iff} $a$ is identical to $b$
	\end{enumerate}
Here, `$\relstroke$' indicates \emph{pairwise} indiscernibility; `$\monadstroke$' indicates \emph{monadic} indiscernibility; and no subscript indicates \emph{complete} indiscernibility. 
\end{define}\noindent
There are several alternative characterisations of $\foi$, two of which will prove useful (see \citeauthor{CasanovasEtAl:EEEFL} \cite[p.\ 508]{CasanovasEtAl:EEEFL} and \citeauthor{Ketland:II} \cite[Theorem 3.17]{Ketland:II}):
	\begin{lem}\label{lem:AlternateCharacterisationOfFOI}
                For any $\langmaybe$-structure $\model{M}$, the following are equivalent:
                \begin{enumerate}[nolistsep]
                		\item $a \foi b$ in $\model{M}$
			\item\label{alternate:atomic} $\model{M} \models \phi(a, \overline{e}) \liff \phi(b, \overline{e})$ for all $n < \omega$, all {atomic} $\phi \in \langno_{n+1}$ and all $\overline{e} \in M^n$
			\item\label{alternate:weak} $\model{M} \models \phi(a, a) \liff \phi(a, b)$ for all $\phi \in \langno_2$ 
		\end{enumerate}
	\end{lem}\noindent 
Quine was the first philosopher to analyse all six grades of indiscernibility. His fullest discussion of them ended as follows:
	\begin{quote}
		May there even be many intermediate grades? The question is ill defined. By imposing special conditions on the form or content of the open sentence used in discriminating two objects, we could define any number of intermediate grades of discriminability, subject even to no linear order. What I have called moderate discriminability [i.e.\ $\yesrel$ or $\norel$], however, is the only intermediate grade that I see how to define at our present high level of generality. \cite[p.\ 116]{Quine:GD}
	\end{quote}
Quine was right that Definition \ref{def:Indiscernibility} essentially exhausts all of the grades of discrimination that are fairly natural, highly general, and which can be defined in terms of satisfaction of $\langno$- and $\langyes$-formulas.\footnote{Though \citeauthor{CaultonButterfield:KILM} \cite{CaultonButterfield:KILM} and \citeauthor{LadymanLinneboPettigrew:IDPL} \cite {LadymanLinneboPettigrew:IDPL} also explore the case of quantifier-free formulas.} Nevertheless, other grades of discrimination are quite natural; we just need to consider alternative methods of definition. (The sense in which they are `intermediate' grades will become clear in \S\ref{s:Entailments} and, as Quine conjectured, we will see that they are not linearly ordered.)

In particular, I shall introduce grades of discrimination that are defined in terms of \emph{isomorphisms}. As a reminder:
\begin{define}
	Let $\model{M}, \model{N}$ be $\langmaybe$-structures. An \emph{isomorphism} from $\model{M}$ to $\model{N}$ is any bijection $\pi : M \longrightarrow N$ such that:
	\begin{enumerate}[nolistsep]
		\item $\overline{e} \in R^\model{M}$ iff $\pi(\overline{e}) \in R^\model{N}$, for all $n$-place $\langmaybe$-predicates $R$ and all $\overline{e} \in M^n$
		\item $\pi(f^\model{M}(\overline{e})) = f^\model{N}(\pi(\overline{e}))$, for all $n$-place $\langmaybe$-function-symbols $f$ and all $\overline{e} \in M^n$
	\end{enumerate}
A \emph{symmetry} on $\model{M}$ is an isomorphism from $\model{M}$ to $\model{M}$. 
\end{define}\noindent
Isomorphisms preserve $\langyes$-formulas (see e.g.\ \citeauthor{Marker:MT} \cite[pp.\ 13--14]{Marker:MT}):
\begin{lem}\label{lem:IsomorphismPreservation}
		Let $\model{M}, \model{N}$ be $\langmaybe$-structures, and $\pi : \model{M} \functionto \model{N}$ be an isomorphism. For all $n < \omega$, all $\phi \in \langyes_n$ and all $\overline{e} \in M^n$: 
		$$\model{M} \models \phi(\overline{e})\text{ iff }\model{N}\models \phi(\pi(\overline{e}))$$
\end{lem}\noindent
There is therefore a good sense in which objects linked by a symmetry cannot be discriminated. Consequently, symmetries are a source of grades of discrimination, and I shall be interested in three distinct \emph{grades of symmetry}:
\begin{define}
	For any $\langmaybe$-structure $\model{M}$ with $a, b \in M$:
	\begin{enumerate}[nolistsep]
		\item $a \symbare b$ in $\model{M}$ iff there is a symmetry $\pi$ on $\model{M}$ with $\pi(a) = b$
		\item $a \sympair b$ in $\model{M}$ {iff} there is a symmetry $\pi$ on $\model{M}$ with $\pi(a) = b$ and $\pi(b) = a$
		\item $a \symtotal b$ in $\model{M}$ {iff} there is a symmetry $\pi$ on $\model{M}$ with $\pi(a)= b$, $\pi(b) = a$ and $\pi(x) = x$ for all $x \notin \{a, b\}$	
	\end{enumerate}
\end{define}\noindent
These three grades have already received some philosophical attention;\footnote{$\symbare$ is considered by \citeauthor{Ketland:II} \cite{Ketland:SII}, \cite{Ketland:II} under the name `structural indiscernibility', and by \citeauthor{LadymanLinneboPettigrew:IDPL} \cite{LadymanLinneboPettigrew:IDPL} under the name `symmetry'.  $\sympair$ is considered by \citeauthor{LadymanLinneboPettigrew:IDPL} \cite{LadymanLinneboPettigrew:IDPL} under the name  `full symmetry'. $\symtotal$ is considered by \citeauthor{Ketland:II} \cite{Ketland:SII}, \cite{Ketland:II}, who writes `$\pi_{ab}$' to indicate that $a \symtotal b$.} one of my aims is to incorporate them into the discussion in a systematic way.

In defining the notion of an isomorphism, the only object-language symbols which are mentioned are those of the signature; there is no need to mention `$=$'. Nevertheless, the notion of an isomorphism---and hence each grade of symmetry---straightforwardly depends upon the notion of identity. After all, an isomorphism is a bijection, which is to say it maps \emph{unique} objects to \emph{unique} objects, and \emph{vice versa}. This dependence on identity is reflected in Lemma \ref{lem:IsomorphismPreservation}: symmetries preserve $\langyes$-formulas. 

If we want to avoid treating identity as a primitive---for philosophical or technical reasons---then the notion of an isomorphism is probably too strong. In looking for a weaker notion, a first thought would be to consider functions between structures that need not be bijections. (In this regard, \emph{strict homomorphisms} are sometimes considered.) But this is insufficiently concessive, since the very idea of a \emph{function} presupposes the notion of identity, for a function maps each object (or $n$-tuple) to a \emph{unique} object. Instead, then, we should consider structure-preserving \emph{relations} that may hold between structures. The appropriate notion is provided by \citeauthor{CasanovasEtAl:EEEFL} \cite[Definition 2.5]{CasanovasEtAl:EEEFL}; recall from \S\ref{s:Preliminaries} that $\overline{d}\Pi\overline{e}$ abbreviates $\langle d_{1}, e_{1}\rangle, \ldots, \langle d_{n}, e_{n}\rangle \in \Pi$:\footnote{They credit a special case of this to \citeauthor{BlokPigozzi:PL} \cite[p.\ 343]{BlokPigozzi:PL}; see also \cite[\S5]{Waszkiewicz:NII}.}
\begin{define}
	Let $\model{M}, \model{N}$ be $\langmaybe$-structures. A \emph{relativeness correspondence from $\model{M}$ to $\model{N}$} is any relation $\Pi \subseteq M \times N$ with $\domain(\Pi) = M$ and $\range(\Pi) = N$ such that:
	\begin{enumerate}[nolistsep]
		\item $\overline{d} \in R^\model{M}$ {iff} $\overline{e} \in R^\model{N}$, for all $n$-place $\langmaybe$-predicates $R$ and all $\overline{d}\Pi\overline{e}$
		\item $f^\model{M}(\overline{d}) \Pi f^\model{N}(\overline{e})$, for all $n$-place $\langmaybe$-function-symbols $f$ and all $\overline{d}\Pi\overline{e}$
	\end{enumerate}
A \emph{relativity} on $\model{M}$ is a relativeness correspondence from $\model{M}$ to $\model{M}$. 
\end{define}\noindent
\citeauthor{CasanovasEtAl:EEEFL} \cite[Proposition 2.6]{CasanovasEtAl:EEEFL} show that relativeness correspondences preserve $\langno$-formulas:
	\begin{lem}\label{lem:RelativenessPreservation}
		Let $\model{M}, \model{N}$ be $\langmaybe$-structures, and $\Pi$ be a relativeness correspondence from $\model{M}$ to $\model{N}$. For all $n < \omega$, all $\phi \in \langno_n$, and all $\overline{d}\Pi\overline{e}$: 
		$$\model{M} \models \phi(\overline{d})\text{ iff }\model{N}\models \phi(\overline{e}) $$
\end{lem}\noindent
There is therefore a good sense in which objects linked by a relativity cannot be discriminated. So, by simple analogy with the three grades of symmetry, I shall consider three \emph{grades of relativity}:
\begin{define}
	For any $\langmaybe$-structure $\model{M}$ with $a, b \in M$:
	\begin{enumerate}[nolistsep]
		\item $a \relativebare b$ in $\model{M}$ {iff} there is a relativity $\Pi$ on $\model{M}$ with $a \Pi b$
		\item $a \relativepair b$  in $\model{M}$ {iff} there is a relativity $\Pi$ on $\model{M}$ with $a \Pi b$ and $b \Pi a$
		\item $a \relativetotal b$ in $\model{M}$  {iff} there is a relativity $\Pi$ on $\model{M}$ with $a \Pi b$, $b \Pi a$ and $x \Pi x$ for all $x \nfoi a$ and $x \nfoi b$ 
	\end{enumerate}
\end{define}\noindent
Unlike the grades of symmetry, the grades of relativity have not yet been considered by philosophers interested in grades of discrimination. However, there is no principled reason for this omission. Indeed, a central claim of this paper is that relativeness correspondences (and hence grades of relativity) are the appropriate $\langno$-surrogate for isomorphisms (and hence grades of symmetry). This claim should already be plausible, given that we arrived at the notion of a relativeness correspondence by relaxing the notion of an isomorphism, and given the immediate comparison between Lemmas \ref{lem:IsomorphismPreservation} and \ref{lem:RelativenessPreservation}. The claim will receive further support during this paper.

For the reader's convenience, the following table summarises the twelve grades of discrimination: 

\begin{table}[h]
\setlength\tabcolsep{0.025\textwidth}
\begin{tabular}{@{}p{0.05\textwidth}p{0.4\textwidth}p{0.45\textwidth}@{}}
			Grade & Informal gloss & Definition sketch \\
			\hline
			 $=$ & genuine identity & $a = b$\\
			$\yesrel$ & pairwise $\langyes$-indiscernibility & $\phi(a,b) \liff \phi(b,a)$, all $\phi\in\langyes_2$\\
			$\yesmonad$ & monadic $\langyes$-indiscernibility & $\phi(a) \liff \phi(b)$, all $\phi \in \langyes_1$  \\\\
			
			$\foi$ & complete $\langno$-indiscernibility & $\phi(a, \overline{e}) \liff \phi(b, \overline{e})$, all $\overline{e}$ and $\phi\in \langno_n$\\
			 $\norel$ & pairwise $\langno$-indiscernibility & $\phi(a,b) \liff \phi(b,a)$, all $\phi\in\langno_2$\\
			$\nomonad$ & monadic $\langno$-indiscernibility & $\phi(a) \liff \phi(b)$, all $\phi \in \langno_1$  \\\\
			
			 $\symtotal$& complete symmetry & a permutation $\pi(a) = b$, $\pi(b) = a$ and $\pi(x) = x$ all $x \notin \{a, b\}$\\
			$\sympair$& pairwise symmetry & a permutation $\pi(a) = b$, $\pi(b) = a$\\
			$\symbare$& monadic symmetry & a permutation $\pi(a)= b$\\\\
			
			$\relativetotal$ & complete relativity & a relativity $a\Pi  b$, $b\Pi a$ and $x \Pi x$ all $x \nfoi a, x \nfoi b$\\
			$\relativepair$& pairwise relativity & a relativity $a \Pi b$, $b \Pi a$ \\
			$\relativebare$& monadic relativity & a relativity $a \Pi b$ 
\end{tabular}
\end{table}
 
\section{Entailments between the grades}\label{s:Entailments}
Having defined twelve grades of discrimination, my first task is to characterise the relationships between them. More precisely: I shall build upon some existing results (mentioned in endnotes) to provide a complete account of the entailments and non-entailments between the various grades of discrimination. 

For any two grades of discrimination $\mathrm{R}$ and $\mathrm{S}$, say that $\mathrm{R}$ \emph{entails} $\mathrm{S}$ \emph{iff} for any structure $\model{M}$ and any $a, b \in M$, if $a \mathrm{R} b$ then $a \mathrm{S} b$ in $\model{M}$. Entailment is relativised to particular classes of structures---e.g.\ to structures with relational signatures---in the obvious way. In \S\ref{s:FinitaryCase}, I shall consider the special case of entailments where we restrict our attention to \emph{finite} structures. However, the target result for this section is the general case:
\addtocounter{dummy}{2}
	\begin{thm}[Entailments between the grades]\label{thm:BigMap} These Hasse Diagrams characterise the entailments between our grades of discrimination:
	
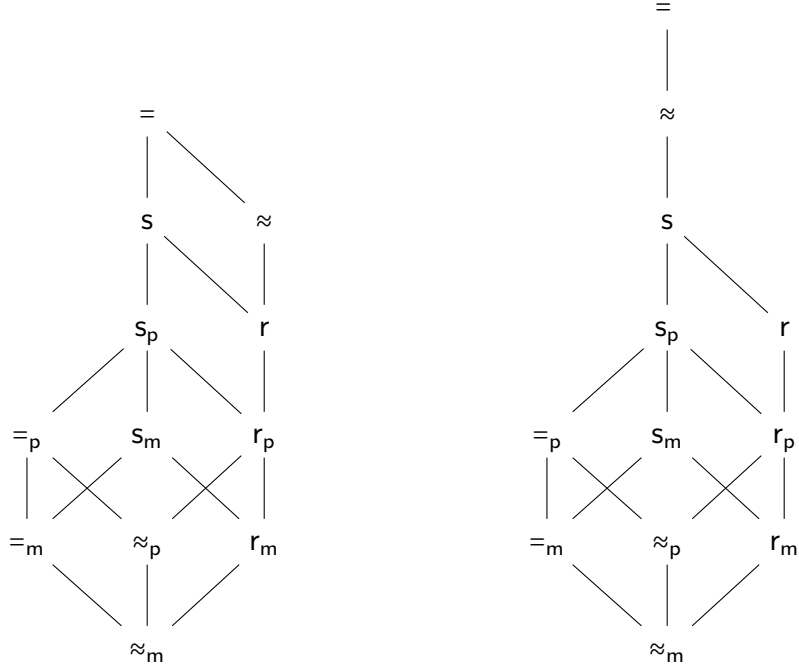
\begin{quote}
\begin{tikzpicture}
		\matrix (m) [matrix of math nodes, row sep=2em, 
		column sep=2em, text height=1.5ex, text depth=0.25ex] 
		{  &  & \\
		   &= & \\
		 & \symtotal & \foi  \\
		 & \sympair & \relativetotal \\
		 \yesrel & \symbare & \relativepair \\
		  \yesmonad & \norel & \relativebare \\
		 & \nomonad&  \\};
		\draw(m-2-2)--(m-3-2)--(m-4-2)--(m-5-2)--(m-6-1)--(m-7-2);
		\draw(m-2-2)--(m-3-3)--(m-4-3)--(m-5-3)--(m-6-3)--(m-7-2);
		\draw(m-3-2)--(m-4-3);
		\draw(m-4-2)--(m-5-3);
		\draw(m-5-2)--(m-6-3);
		\draw(m-4-2)--(m-5-1)--(m-6-2)--(m-7-2);
		\draw(m-5-1)--(m-6-1);
		\draw(m-5-3)--(m-6-2);
		\end{tikzpicture}
\hfill
\begin{tikzpicture}
		\matrix (m) [matrix of math nodes, row sep=2em, 
		column sep=2em, text height=1.5ex, text depth=0.25ex] 
		{  & = & & \\
		   & \foi & & \\
		 & \symtotal & &  \\
		 & \sympair & \relativetotal& \\
		 \yesrel & \symbare & \relativepair \\
		  \yesmonad & \norel & \relativebare \\
		 & \nomonad&  \\};
		\draw(m-1-2)--(m-2-2)--(m-3-2)--(m-4-2)--(m-5-2)--(m-6-1)--(m-7-2);
		\draw(m-4-3)--(m-5-3)--(m-6-3)--(m-7-2);
		\draw(m-3-2)--(m-4-3);
		\draw(m-4-2)--(m-5-3);
		\draw(m-5-2)--(m-6-3);
		\draw(m-4-2)--(m-5-1)--(m-6-2)--(m-7-2);
		\draw(m-5-1)--(m-6-1);
		\draw(m-5-3)--(m-6-2);
		\end{tikzpicture}
	\end{quote}
	The left diagram considers the case of arbitrary signatures; the right diagram considers entailment when restricted to structures with relational signatures.
\end{thm}\addtocounter{dummy}{-3}\noindent
To explain the notation: there is a path down the page from $\text{R}$ to $\text{S}$ \emph{iff} $\text{R}$ entails $\text{S}$. So in the case of arbitrary signatures: $=$ entails $\foi$; $\foi$ does not entail $=$; $\symtotal$ does not entail $\foi$; and $\foi$ does not entail $\symtotal$. In the case of relational signatures: $\foi$ entails $\symtotal$; hence $\foi$ entails $\sympair$; etc. I shall start by proving the entailments:\footnote{\citeauthor{Quine:GD} \cite{Quine:GD} proves case (\ref{yesyesnonohierarchy}); see also \citeauthor{Ketland:SII} \cite[p.\ 307]{Ketland:SII}, \cite[\S3.2]{Ketland:II}, \citeauthor{LadymanLinneboPettigrew:IDPL} \cite[Theorem 5.2]{LadymanLinneboPettigrew:IDPL}. \citeauthor{CaultonButterfield:KILM} \cite[Theorem 1]{CaultonButterfield:KILM} prove case (\ref{symtoyes}) when restricted to relational signatures; see also \citeauthor{Ketland:II} \cite[Lemma 3.22]{Ketland:II} and \citeauthor{LadymanLinneboPettigrew:IDPL} \cite[Theorem 9.17, 9.20]{LadymanLinneboPettigrew:IDPL}. \citeauthor{Ketland:II} \cite[Theorem 3.23]{Ketland:II} proves case (\ref{foitosymtotal}).}
\begin{lem}\label{lem:BigMapEntailments}For structures with arbitrary signatures:
	\begin{enumerate}[nolistsep]
		\item\label{identityattop} $=$ entails both $\foi$ and $\symtotal$
		\item\label{foitonorel} $\foi$ entails $\relativetotal$
		\item\label{yesnoyesonhierarchy} $\yesrel$ entails $\norel$, and $\yesmonad$ entails $\nomonad$
		\item\label{yesyesnonohierarchy} $\yesrel$ entails $\yesmonad$, and $\norel$ entails $\nomonad$

		\item\label{symtoyes} $\sympair$ entails $\yesrel$, and $\symbare$ entails $\yesmonad$	
		\item\label{reltono} $\relativepair$ entails $\norel$, and $\relativebare$ entails $\nomonad$

		\item\label{StrengthSymmetryRelativity:total} $\symtotal$ entails $\relativetotal$, $\sympair$ entails $\relativepair$, and $\symbare$ entails $\relativebare$

		\item\label{symhierarchy} $\symtotal$ entails $\sympair$, and $\sympair$ entails $\symbare$
		\item\label{relativehierarchy} $\relativetotal$ entails $\relativepair$, and $\relativepair$ entails $\relativebare$
	\end{enumerate}
	For structures with relational signatures, but not in general:  
	\begin{enumerate}[nolistsep]\addtocounter{enumi}{9}
		\item \label{foitosymtotal} $\foi$ entails $\symtotal$
	\end{enumerate}
	\begin{proof}
		(\ref{identityattop}). The identity map is a symmetry.
		
		(\ref{foitonorel}). The relation given by $x \Pi x$ iff $x \foi x$ is a relativity.
		
		(\ref{yesnoyesonhierarchy})--(\ref{yesyesnonohierarchy}). Immediate from the definitions.

		(\ref{symtoyes})--(\ref{reltono}). Immediate from Lemmas \ref{lem:IsomorphismPreservation} and \ref{lem:RelativenessPreservation}.
		
		(\ref{StrengthSymmetryRelativity:total}). Every symmetry can be regarded as a relativity.

		(\ref{symhierarchy})--(\ref{relativehierarchy}). Immediate from the definitions.
		
		(\ref{foitosymtotal}). If $a \foi b$, then for any $n < \omega$, any atomic formula $\phi \in \langno_n$, and any $\overline{e} \in (M\setminus \{a, b\})^n$, we have $\model{M} \models \phi(a, \overline{e}) \liff \phi(b, \overline{e})$. 
	\end{proof}
\end{lem}\noindent
It remains to demonstrate the non-entailments:\footnote{\citeauthor{Ketland:II} \cite[p.\ 8]{Ketland:II} and \citeauthor{LadymanLinneboPettigrew:IDPL} \cite[Theorem 5.2]{LadymanLinneboPettigrew:IDPL} use $\model{\Blank}$ to prove case (\ref{foinotidentity}); see also \citeauthor{Button:RSIC} \cite[p.\ 218]{Button:RSIC} and \citeauthor{Ketland:SII} \cite[p.\ 309]{Ketland:SII}. \citeauthor{LadymanLinneboPettigrew:IDPL} \cite[Theorem 9.17]{LadymanLinneboPettigrew:IDPL} use $\model{\Black}$ to
prove case (\ref{symtotalnotfoi}), noting that it is the analogue of \citeauthor{Black:II}'s \cite*{Black:II} two-sphere world. \citeauthor{Button:RSIC} \cite[p.\ 218]{Button:RSIC}, \citeauthor{Ketland:SII} \cite[p.\ 310]{Ketland:SII} and \citeauthor{LadymanLinneboPettigrew:IDPL} \cite[Theorem 7.12]{LadymanLinneboPettigrew:IDPL} use an example like $\model{\ThreeSphere}$. \citeauthor{LadymanLinneboPettigrew:IDPL} \cite[Theorem 5.2]{LadymanLinneboPettigrew:IDPL} use an example like $\model{\Rotational}$.  \citeauthor{CaultonButterfield:KILM} \cite[pp.\ 60--2]{CaultonButterfield:KILM} and \citeauthor{LadymanLinneboPettigrew:IDPL} \cite[Theorem 9.20]{LadymanLinneboPettigrew:IDPL} use examples like $\model{\TwoBallBig}$.}
\begin{lem}\label{lem:BigMapCounterEntailments}
	For structures with relational signatures:
		\begin{enumerate}[nolistsep]
			\item\label{foinotidentity} $\foi$ does not entail $=$
			\item\label{symtotalnotfoi} $\symtotal$ does not entail $\foi$
			\item\label{relativetotalnotyesmonad} $\relativetotal$ does not entail $\yesmonad$
			\item\label{symbarenotnopair} $\symbare$ does not entail $\norel$
			\item\label{sympairnotrelativetotal} $\sympair$ does not entail $\relativetotal$
			\item\label{yesrelnotrelativebare} $\yesrel$ does not entail $\relativebare$
		\end{enumerate}
	Moreover, for structures with arbitrary signatures:
		\begin{enumerate}[nolistsep]\addtocounter{enumi}{6}
			\item\label{foinotyesmonad} $\foi$ does not entail $\yesmonad$
		\end{enumerate}
	\begin{proof}
		(\ref{foinotidentity}). In this unlabelled graph, $\model{\Blank}$, we have $1 \foi 2$ but $1 \neq 2$:
			\begin{quote}
			\centering
			\begin{tikzpicture}[descr/.style={fill=white,inner sep=2.5pt}] 
			\matrix (m) [matrix of math nodes, row sep=2em, 
			column sep=2em, text height=1.5ex, text depth=0.25ex] 
			{1 & 2\\ }; 
			\end{tikzpicture} 
			\end{quote}
			
		(\ref{symtotalnotfoi}). In this unlabelled graph, $\model{\Black}$, we have $1 \symtotal 2$ but $1 \nfoi 2$:
			\begin{quote}
			\centering
			\begin{tikzpicture}[descr/.style={fill=white,inner sep=2.5pt}] 
			\matrix (m) [matrix of math nodes, row sep=2em, 
			column sep=2em, text height=1.5ex, text depth=0.25ex] 
			{1 & 2\\ }; 
			\path[-, thick] (m-1-1) edge (m-1-2);
			\end{tikzpicture} 
			\end{quote}

		(\ref{relativetotalnotyesmonad}). In this unlabelled graph, $\model{\ThreeSphere}$, we have $1 \relativetotal 2$ but $1 \nyesmonad 2$:
			\begin{quote}
			\centering
			\begin{tikzpicture}[descr/.style={fill=white,inner sep=2.5pt}] 
			\matrix (m) [matrix of math nodes, row sep=2em, 
			column sep=2em, text height=1.5ex, text depth=0.25ex] 
			{\phantom{3}& 1 & 2 & 3\\ }; 
			\path[-, thick] (m-1-3) edge (m-1-2);
			\path[-, thick] (m-1-3) edge (m-1-4);
			\end{tikzpicture} 
			\end{quote}
			
		(\ref{symbarenotnopair}). In this unlabelled directed graph, $\model{\Rotational}$, we have $1 \symbare 2$ but $1 \nnorel 2$:
			
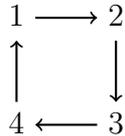
\begin{quote}
			\centering
			\begin{tikzpicture}[descr/.style={fill=white,inner sep=2.5pt}] 
			\matrix (m) [matrix of math nodes, row sep=2em, 
			column sep=2em, text height=1.5ex, text depth=0.25ex] 
			{1 & 2\\ 
			 4 & 3\\ }; 
			\path[->, thick] (m-1-1) edge (m-1-2);
			\path[->, thick] (m-1-2) edge (m-2-2);
			\path[->, thick] (m-2-2) edge (m-2-1);
			\path[->, thick] (m-2-1) edge (m-1-1);
			\end{tikzpicture}
			\end{quote}
					
		(\ref{sympairnotrelativetotal}). In $\model{\Rotational}$, again, we have $1 \sympair 3$ but $1 \nrelativetotal 3$.
		
		(\ref{yesrelnotrelativebare}). Let $\model{\TwoBallBig}$ be the disjoint union of a complete countably-infinite graph with a complete uncountable graph, i.e.:
		\begin{align*}
			\TwoBallBig &:= \mathbb{R}\\
			R^{\model{\TwoBallBig}} &:= \{\langle n, m\rangle \in \mathbb{N}^2 \mid n \neq m\} \cup \{\langle p, q \rangle \in (\mathbb{R} \setminus \mathbb{N})^2 \mid p \neq q\}
		\end{align*}
		By taking a Skolem Hull of $\model{\TwoBallBig}$ containing $1\in \mathbb{N}$ and any $e \in \mathbb{R}\setminus\mathbb{N}$, we see that $1 \yesrel e$. Now suppose that $\Pi$ is a relativity with $1 \Pi e$. Since $\Pi$ must preserve the edges of the graph, and every element in either `cluster' has an edge to every element in the cluster except itself, $\Pi$ must be a bijection between $\mathbb{N}$ and $\mathbb{R}\setminus \mathbb{N}$. Contradiction; so $1 \nrelativebare e$.
		
		(\ref{foinotyesmonad}). Augment $\model{\Blank}$ by adding a single constant which picks out $1$.
	\end{proof}
\end{lem}\noindent
It is simple to check that Lemmas \ref{lem:BigMapEntailments} and \ref{lem:BigMapCounterEntailments} yield Theorem \ref{thm:BigMap}. This Theorem allows us to compare the consequences of imposing various grades of discrimination as criteria of identity.

I should comment briefly on the philosophical significance of the constructions used in Lemma \ref{lem:BigMapCounterEntailments}. The existence of $\model{\Blank}$ is guaranteed by absolutely standard model theory. However, $\model{\Blank}$ contains two distinct objects that are `blank': from the perspective of $\model{\Blank}$, these objects have no properties or relations to anything, so that their distinctness must be \emph{brute}. And this might suggest that the use of absolutely standard model theory begs the question against anyone who believes in a non-trivial criterion of identity. Fortunately it does not, but it is worth carefully explaining why. 

Let Fran be a philosopher who advocates a non-trivial criterion of identity: in particular, Fran thinks that $x$ and $y$ are identical iff $x \foi y$. However, bearing in mind the discussion of \S\ref{s:Preliminaries}---particularly of a signature which allows us only to describe eye colour---Fran advances this criterion of identity with respect to some \emph{particular} signature, $\lang{F}$. Now, if $\model{\Blank}$ is presented as an $\lang{F}$-structure, then Fran will certainly deny that $\model{A}$ could exist. However, Fran can make sense of $\model{\Blank}$ by regarding it as a $\lang{G}$-structure, where $\lang{G}$ is a signature which is  impoverished compared with $\lang{F}$. Construed thus, $\model{\Blank}$ begs no question against Fran, because it poses no threat to her proposed criterion of identity. 

To be clear: I am not trying to endorse or defend Fran's position.\footnote{I was once on Fran's side \cite[p.\ 220]{Button:RSIC}; but I have changed my mind \cite[p.\ 211 n.\ 8]{Button:LoR}.} My point is simply that everyone, including Fran, can make sense of standard model theory.

\section{A Galois Connection}\label{s:Galois}
Theorem \ref{thm:BigMap} graphically demonstrates that grades of symmetry are to grades of $\langyes$-indiscernibility as grades of relativity are to grades of $\langno$-indiscernibility. In this section, I develop this point by outlining a Galois Connection between isomorphisms and relativeness correspondences. (The results of this section can be fruitfully compared with those of \citeauthor{BonnayEngstrom:ID} \cite{BonnayEngstrom:ID}; we discovered our results independently.)

Lemma \ref{lem:IsomorphismPreservation} has an obvious converse: every bijective map which preserves all $\langyes$-formulas is an isomorphism. However, there is no converse to Lemma \ref{lem:RelativenessPreservation}. To make this more precise, consider the following definition:
\begin{define}\label{def:NearCorrespondence}
	Let $\model{M}, \model{N}$ be $\langmaybe$-structures. A \emph{near-correspondence from $\model{M}$ to $\model{N}$} is any relation $\Pi \subseteq M \times N$ with $\domain(\Pi) = M$ and $\range(\Pi) = N$ such that, for all $n < \omega$, all $\phi \in \langno_n$, and all $\overline{d}\Pi\overline{e}$:
		\begin{align*}
			\model{M} \models \phi(\overline{d})&\text{ iff }\model{N} \models \phi(\overline{e})
		\end{align*}
\end{define}\noindent
Lemma \ref{lem:RelativenessPreservation} states that every relativeness correspondence is a near-correspondence. But the converse fails. Let $\model{\NearCorrespondence}$ be an $\{f\}$-structure, defined as follows:
	\begin{align*}
		F &= \{1, 2\}\\
		f^{\model{\NearCorrespondence}}(1) &= f^{\model{\NearCorrespondence}}(2) = 2
	\end{align*}
Then $\Pi = \{\langle 1, 2\rangle, \langle 2, 1\rangle\}$ is clearly a near-correspondence from $\model{\NearCorrespondence}$ to $\model{\NearCorrespondence}$, but not a relativeness correspondence.

However, there is an elegant connection between near-correspondences (and hence relativeness correspondences) and isomorphisms on the models we obtain by quotienting using $\foi$.  The use of such quotients is standard in model theory without identity, and the central idea is summed up in the following Definition and Lemma (see \citeauthor{CasanovasEtAl:EEEFL} \cite[Definition 2.3--2.4]{CasanovasEtAl:EEEFL}):\footnote{\citeauthor{CasanovasEtAl:EEEFL} trace the definition and lemma back to \citeauthor{Monk:ML} \cite[presumably Exercises 29.33--34]{Monk:ML}. This has recently been rediscovered by philosophers, e.g.\ Ketland \cite[p.\ 307 n.\ 10]{Ketland:SII}, \cite[Theorem 3.12]{Ketland:II}.} 
\begin{define}
	Let $\model{M}$ be any $\langmaybe$-structure. Then $\quotient{\model{M}}$ is the $\langmaybe$-structure obtained by quotienting $\model{M}$ by $\foi$. We denote its members with $\quotient{a}_\model{M} = \{b \in M \mid a \foi b \text{ in }\model{M}\}$ and, when no confusion can arise, we dispense with the subscript, talking of $\quotient{a}$ rather than $\quotient{a}_\model{M}$. Now $\quotient{\model{M}}$ is defined as follows:
	\begin{align*}
		\quotient{M} & = \{\quotient{a} \mid a \in M\}\\
		R^{\quotient{\model{M}}} & = \{\quotient{\overline{e}} \in \quotient{M}^n \mid \overline{e} \in R^\model{M}\} && \text{all }n\text{-place }\langmaybe\text{-predicates }R\\
		f^{\quotient{\model{M}}}(\quotient{\overline{e}}) & = \quotient{f^{\model{M}}(\overline{e})} && \text{all }n\text{-place }\langmaybe\text{-function-symbols }f\text{ and all }\overline{e} \in M^n
	\end{align*}
	\end{define}
\begin{lem}\label{lem:Quotient}Let $\model{M}$ be an $\langmaybe$-structure. For all $n < \omega$, all $\phi \in \langno_n$ and all $\overline{e} \in M^n$:
	\begin{align*}
		\model{M} \models \phi(\overline{e}) & \text{ iff }\quotient{\model{M}} \models \phi(\quotient{\overline{e}})
	\end{align*}	
\end{lem}\noindent
\citeauthor{CasanovasEtAl:EEEFL} \cite[Proposition 2.6]{CasanovasEtAl:EEEFL} note that there is a relativeness correspondence from $\model{M}$ to $\model{N}$ iff $\quotient{\model{M}} \isomorphic \quotient{\model{N}}$. I wish to build on this; and I begin with some definitions:
\begin{define}\label{define:Galois}
	For any $\langmaybe$-structures $\model{M}, \model{N}$:
		\begin{enumerate}[nolistsep]
			\item $\GaloisNear(\model{M}, \model{N})$ is the set of near-correspondences from $\model{M}$ to $\model{N}$
			\item  $\GaloisIso(\model{M}, \model{N})$ is the set of isomorphisms from $\quotient{\model{M}}$ to $\quotient{\model{N}}$
			\item $\galup(\model{M}, \model{N}) : \GaloisIso(\model{M}, \model{N}) \longrightarrow \GaloisNear(\model{M}, \model{N})$ is given by: $a \pi^\galup b$ iff $\pi(\quotient{a}) = \quotient{b}$
	 		\item $\galdown(\model{M}, \model{N}) : \GaloisNear(\model{M}, \model{N}) \longrightarrow \GaloisIso(\model{M}, \model{N})$ is given by: $\Pi^\galdown(\quotient{a}) = \quotient{b}$ iff there are $a' \foi a$ and $b' \foi b$ such that $a' \Pi b'$
	\end{enumerate}
	Say that $\Pi \in \GaloisNear(\model{M}, \model{N})$ is \emph{maximal} iff no strict superset of $\Pi$ is in $\GaloisNear(\model{M}, \model{N})$.
\end{define}\noindent
I prove that these are genuine definitions, i.e.\ that $\galup(\model{M}, \model{N})$ and $\galdown(\model{M}, \model{N})$ are functions. I begin with $\galup(\model{M}, \model{N})$:
\begin{lem}\label{lem:QuotientSymmetryGivesRelativity}
	If $\pi \in \GaloisIso(\model{M}, \model{N})$, then $\pi^\galup$ is a relativeness correspondence, and hence $\pi^{\galup} \in \GaloisNear(\model{M}, \model{N})$. 
	\begin{proof}
		Fix $n < \omega$ and suppose that $\overline{d}\pi^\galup\overline{e}$; i.e.\ that $\pi(\overline{\quotient{d}}) = \overline{\quotient{e}}$. For each $n$-place $L$-predicate $R$, observe:
			\begin{align*}
				\overline{d} \in R^\model{M}
					\text{ iff }\overline{\quotient{d}} \in R^{\quotient{\model{M}}}
					\text{ iff }\pi(\overline{\quotient{d}}) \in R^{\quotient{\model{N}}}
					\text{ iff }\overline{\quotient{e}} \in R^{\quotient{\model{N}}}
					\text{ iff }\quotient{e} \in R^\model{N}
			\end{align*}
		For each $n$-place $L$-function-symbol $f$, observe:
			\begin{align*}
				\pi(\quotient{f^\model{M}(\overline{d})}) =
				\pi(f^{\quotient{\model{M}}}(\overline{\quotient{d}})) =
				f^{\quotient{\model{N}}}(\pi(\overline{\quotient{d}})) = 
				f^{\quotient{\model{N}}}(\overline{\quotient{e}}) = 
				\quotient{f^\model{N}(\overline{e})}
			\end{align*}
		so that $f^\model{M}(\overline{d})\pi^\galup f^\model{N}(\overline{e})$. Hence $\pi^\galup$ is a relativeness correspondence, and so a near-correspondence by Lemma \ref{lem:RelativenessPreservation}.
	\end{proof}
\end{lem}\noindent
To show that $\galdown(\model{M}, \model{N})$ is a function, we need a subsidiary result:
\begin{lem}\label{lem:RelativitiesAndFoi}
	Let $\Pi$ be a near-correspondence from $\model{M}$ to $\model{N}$, with $a \Pi b$ and $a' \Pi b'$. Then $a \foi a'$ in $\model{M}$ iff $b \foi b'$ in $\model{N}$. 
	\begin{proof}
		Suppose $a \foi a'$ in $\model{M}$; then using Lemmas \ref{lem:RelativenessPreservation} and \ref{lem:AlternateCharacterisationOfFOI}:
			\begin{align*}
				\model{N} \models \phi(b,b) \text{ iff }\model{M} \models \phi(a, a)\text{ iff }\model{M} \models \phi(a, a') \text{ iff }\model{N} \models \phi(b,b')
			\end{align*}
		So $b \foi b'$ in $\model{N}$, by Lemma \ref{lem:AlternateCharacterisationOfFOI}. The converse is similar.
	\end{proof} 
\end{lem}\noindent
It follows that $\galdown(\model{M}, \model{N})$ is a function:
\begin{lem}\label{lem:GalDownGivesIsomorphism}
If $\Pi \in \GaloisNear(\model{M}, \model{N})$, then $\Pi^{\galdown} \in \GaloisIso(\model{M}, \model{N})$.
	\begin{proof} Lemma \ref{lem:RelativitiesAndFoi} immediately yields that $\Pi^\galdown(\overline{a})$ is a well-defined function, and indeed an injection. $\Pi^{\galdown}$ is a surjection, since $\range(\Pi) = N$. It remains to show that $\Pi^\galdown$ preserves structure. For the remainder of the proof, fix $n < \omega$ and ${\overline{a}}, {\overline{b}} \in {M}^n$ such that $\Pi^\galdown(\quotient{\overline{a}}) = \quotient{\overline{b}}$.

	Let $R$ be any $n$-place $\langmaybe$-predicate.  For each $1 \leq i \leq n$ we have $a'_i \foi {a_i}$ and $b'_i \foi {b_i}$ such that $a'_i \Pi b'_i$; and hence:
				\begin{align*}
					\quotient{\overline{a}} \in R^{\quotient{\model{M}}} &\text{ iff }\overline{a'} \in R^\model{M} \text{ iff }\overline{b'} \in R^\model{N}\text{ iff }\quotient{\overline{b}} \in R^{\quotient{\model{N}}} \text{ iff }\Pi^\galdown(\quotient{\overline{a}})  \in R^{\quotient{\model{N}}} 
				\end{align*}
	Let $f$ be any $n$-place $\langmaybe$-function-symbol. For all $\phi(v_1, v_2) \in \langno_{2}$, define $\phi_{f}(v_1, \overline{x})$ as $\phi(v_1, f(\overline{x}))$; and let $b' \in N$ be such that $f^\model{M}(\overline{a})\Pi b'$. Invoking Lemma \ref{lem:AlternateCharacterisationOfFOI}:
				\begin{align*}
					\model{N} \models \phi(b', b') &\text{ iff }\model{M} \models \phi(f^\model{M}(\overline{a}), f^\model{M}(\overline{a}))\\
					&\text{ iff }\model{M} \models \phi_f(f^\model{M}(\overline{a}), \overline{a})\\
					&\text{ iff }\model{N} \models \phi_f(b', \overline{b})\\
					&\text{ iff }\model{N} \models \phi(b', f^\model{N}(\overline{b}))
				\end{align*}
			Hence  $b' \foi f^\model{N}(\overline{b})$ by Lemma \ref{lem:RelativenessPreservation}. Now:
		\begin{align*}
			\Pi^\galdown(f^{\quotient{\model{M}}}(\overline{\quotient{a}})) & = \Pi^\galdown(\quotient{f^\model{M}(\overline{a})}) 
			= \quotient {b'}
			= \quotient{f^\model{N}(\overline{b})}
			= f^{\quotient{\model{N}}}(\overline{\quotient{b}})
			= f^{\quotient{\model{N}}}(\Pi^\galdown(\overline{\quotient{a}}))
		\end{align*}
	so that functions are preserved.
	\end{proof}
\end{lem}\noindent
Lemmas \ref{lem:QuotientSymmetryGivesRelativity} and \ref{lem:GalDownGivesIsomorphism} together show that Definition \ref{define:Galois} is a proper definition. Its significance resides in the following:
\begin{thm}[Galois Connection on $\foi$-quotients]\label{thm:GaloisConnection}
	For each $\Pi \in \GaloisNear(\model{M}, \model{N})$ and each $\pi \in \GaloisIso(\model{M}, \model{N})$: $\Pi^\galdown = \pi$ iff $\Pi \subseteq \pi^\galup$.
	\begin{proof}
	\emph{Left-to-right.} Suppose $\Pi^\galdown = \pi$. Fix $\langle d, e \rangle \in \Pi$; then $\pi(\quotient{d}) = \quotient{e}$, so $d \pi^\galup e$. 

	\emph{Right-to-left.} Suppose $\Pi \subseteq \pi^\galup$. Where $\Pi^\galdown(\quotient{d}) = \quotient{e}$, there are $d' \foi {d}$ and $e' \foi {e}$ such that $d' \Pi e'$ and hence $d'\pi^\galup e'$; so $\quotient{e} = \quotient{e'}  = \pi(\quotient{d'}) = \pi(\quotient{d})$. Hence $\Pi^\galdown(\quotient{d}) = \pi(\quotient{d})$, for all $\quotient{d} \in \quotient{M}$. 
	\end{proof}
\end{thm}\noindent
This Theorem highlights the depth of the connection between isomorphisms and relativeness correspondences. Additionally, it strengthens the claim that relativeness correspondences are the $\langno$-analogue of isomorphisms. For, given that there are near-correspondences that are not relativeness correspondences, one might have worried that relativeness correspondences \emph{compete} with the near-correspondences to be the $\langno$-analogue of isomorphism. However, the appearance of competition vanishes, once we consider some consequences of the Galois Connection:
\begin{lem}\label{lem:GaloisConsequences}
	For any $\langmaybe$-structures $\model{M}, \model{N}$:
	\begin{enumerate}[nolistsep]
		\item\label{GaloisConsequences:identity} $\galdown(\model{M}, \model{N}) \circ \galup(\model{M}, \model{N})$ is the identity function
		\item\label{GaloisConsequences:idempotent} $\galup(\model{M}, \model{N}) \circ \galdown(\model{M}, \model{N})$ is idempotent
		\item\label{GaloisConsequences:maximal} If $\pi \in \GaloisIso(\model{M}, \model{N})$, then $\pi^\galup$ is maximal
		\item\label{GaloisConsequences:unique} If $\Pi \in \GaloisNear(\model{M}, \model{N})$, then $(\Pi^\galdown)^\galup$ is the unique maximal relativeness correspondence that extends $\Pi$
	\end{enumerate}
	\begin{proof}
		(\ref{GaloisConsequences:identity})--(\ref{GaloisConsequences:idempotent}). Immediate from the fact that this is a Galois Connection with the partial-ordering on $\GaloisIso(\model{M}, \model{N})$ being identity.
		
		(\ref{GaloisConsequences:maximal}). Let $\Sigma \in \GaloisIso(\model{M}, \model{N})$ be such that $\pi^{\galup} \subseteq \Sigma$, and suppose $a \Sigma b$. By Lemma \ref{lem:GalDownGivesIsomorphism}, $\Sigma^\galdown \in \GaloisNear(\model{M}, \model{N})$, with $\Sigma^\galdown(\quotient{a}) = \quotient{b}$.  For any $d$ such that $a \pi^\galup d$, we have $a \Sigma d$, and hence  $\Sigma^\galdown(\quotient{a}) = \quotient{d}$, so that $b \foi d$. Hence $\pi(\quotient{a}) = \quotient{d} = \quotient{b}$, and so $a \pi^\galup b$. 
			
		(\ref{GaloisConsequences:unique}). Lemma \ref{lem:QuotientSymmetryGivesRelativity}, our Galois Connection, and  case (\ref{GaloisConsequences:maximal}) show that $(\Pi^\galdown)^\galup$ is a maximal relativeness correspondence extending $\Pi$. To show uniqueness, let $\Sigma$ be any maximal relativeness correspondence extending $\Pi$. Consider any $a, b \in M$ such that $a (\Pi^\galdown)^\galup b$. Then there are $a' \foi {a}$, $b' \foi {b}$ such that $a' \Pi b'$, and hence such that $a' \Sigma b'$, since $\Pi \subseteq \Sigma$. Hence, for any $\overline{d}\Sigma\overline{e}$, and any $\phi \in \langno_{n+1}$, by Lemma \ref{lem:RelativenessPreservation}:
			\begin{align*}
				\model{M} \models \phi(a,\overline{d}) &\text{ iff }\model{M} \models \phi(a',\overline{d})\text{ iff }\model{N} \models \phi(b',\overline{e}) \text{ iff }\model{N} \models \phi(b,\overline{e})
			\end{align*}
		Consequently, $\Theta = \Sigma \cup \{\langle a, b\rangle\}$ is a near-correspondence. So $(\Theta^\galdown)^\galup$ is a maximal relativeness correspondence extending $\Sigma$; but $\Sigma$ is itself maximal, so $a \Sigma b$. Generalising, $(\Pi^\galdown)^\galup \subseteq \Sigma$. Since $(\Pi^\galdown)^\galup$ is maximal, $(\Pi^\galdown)^\galup = \Sigma$. 
	\end{proof}
\end{lem}\noindent
The preceding result tells us that every near-correspondence expands to a relativeness correspondence. Accordingly, there is no genuine competition between near-correspondences and relativeness correspondences. Indeed, instead of defining the grades of relativity in terms of relativeness correspondences, we could have defined them in terms of near-correspondences. Or, even more simply, we could have defined them in terms of symmetries on quotient models, as shown by the following  immediate consequence of the preceding results:
\begin{lem}\label{lem:QuotientSymHom} For any $\langmaybe$-structure $\model{M}$:
	\begin{enumerate}[nolistsep]
		\item\label{QuotientSymHom:bare} $a \relativebare b$ in $\model{M}$ \text{iff} $\quotient{a} \symbare \quotient{b}$ in $\quotient{\model{M}}$  
		\item\label{QuotientSymHom:pair} $a \relativepair b$ in $\model{M}$ \text{iff} $\quotient{a} \sympair \quotient{b}$ in $\quotient{\model{M}}$
		\item\label{QuotientSymHom:total} $a \relativetotal b$ in $\model{M}$ \text{iff} $\quotient{a} \symtotal \quotient{b}$ in $\quotient{\model{M}}$
	\end{enumerate}
\end{lem}

\section{Equivalence relations}\label{s:Equivalence}
At various points, I have described the grades of discrimination as behaving like identity. A natural question is whether the grades of discrimination behave like identity in being \emph{equivalence} relations. (Note that I implicitly relied upon the fact that $\foi$ is an equivalence relation in defining the $\foi$-quotient structure.) \citeauthor{LadymanLinneboPettigrew:IDPL} \cite[Theorem 10.22]{LadymanLinneboPettigrew:IDPL} have partially answered this question, in noting that $\yesrel$ and $\norel$ are not transitive (in general). The following result, which employs our Galois Connection, completes the picture.
\begin{thm}\label{lem:EquivalenceIndiscernibility}
	$\yesrel$, $\norel$, $\sympair$ and $\relativepair$ are reflexive and symmetric, but not transitive (in general); the remaining eight grades of discrimination are equivalence relations. 
	\begin{proof}
		Consider the following coloured graph, $\model{\Hexagon}$:
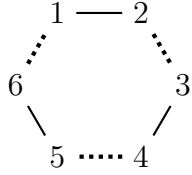
		\begin{quote}\centering
			\begin{tikzpicture}[descr/.style={fill=white,inner sep=2.5pt}] 
			\node (atoma) at (120:1.1){$1$};
			\node (atomb) at (60:1.1){$2$};
			\node (atomc) at (0:1.1){$3$};
			\node (atomd) at (300:1.1){$4$};
			\node (atome) at (240:1.1){$5$};
			\node (atomf) at (180:1.1){$6$};
			\path[-, thick](atoma) edge (atomb);
			\path[-, ultra thick, dotted](atomb) edge (atomc);
			\path[-, thick](atomc) edge (atomd);
			\path[-, ultra thick, dotted](atomd) edge (atome);
			\path[-, thick](atome) edge (atomf);
			\path[-, ultra thick, dotted](atomf) edge (atoma);
			\end{tikzpicture} 
		\end{quote}
		Here, $1 \sympair 2$ and $2 \sympair 3$, but $1 \nnorel 3$. By Theorem \ref{thm:BigMap}, this establishes that $\yesrel$, $\norel$, $\sympair$ and $\relativepair$ are not transitive (in general).
		
		The reflexivity and symmetry of all the grades of indiscernibility are immediate from their definitions, as is the transitivity of $=$, $\foi$, $\yesmonad$ and $\nomonad$.
		
		It is routine to check that all three grades of symmetry are symmetric and reflexive, and that $\symtotal$ and $\symbare$ are transitive. Lemma \ref{lem:QuotientSymHom} entails that the same is true for the respective grades of relativity.
	\end{proof}
\end{thm}\noindent
Since identity is surely transitive, Theorem \ref{lem:EquivalenceIndiscernibility} might seem to provide a knockdown argument against treating any of $\yesrel, \norel, \sympair$ and $\relativepair$ as a criterion of identity. However, this point is a little more subtle than it might initially seem. 

Consider the discussion of $\model{\Blank}$, at the end of \S\ref{s:Entailments}. $\model{\Blank}$ might have seemed to present a counterexample to treating $\foi$ as a criterion of identity. But any philosopher who advocates such a criterion, such as Fran, will maintain that we can make sense of $\model{\Blank}$ by (and only by) treating it as a structure of some artificially \emph{restricted} signature. At that point, $\model{\Blank}$ no longer presents a counterexample to Fran's proposed criterion of identity, which she advances with respect to some \emph{richer} signature.

With this in mind, consider Rach, a philosopher who advocates $\relativepair$ as a criterion of identity. $\model{\Hexagon}$ might seem to pose problems for Rach. But if $\model{\Hexagon}$ is presented with regard to Rach's preferred signature, then it violates her proposed criterion of identity even before we consider issues about transitivity: after all, $\model{\Hexagon}$ is to contain objects which are distinct but (`genuinely') pairwise symmetric. Accordingly, Rach will maintain that we can make sense of $\model{\Hexagon}$ by (and only by) treating it as a structure of some artificially \emph{restricted} signature. And at that point, $\model{\Hexagon}$ no longer demonstrates the non-transitivity of Rach's proposed criterion of identity, which she advances with respect to some \emph{richer} signature.

The situation, then, is slightly odd. From the perspective of anyone who thinks that identity is more fine-grained than any of $\yesrel$, $\norel$, $\sympair$ and $\relativepair$, these four grades of discrimination fail to behave like identity in an absolutely crucial sense, in failing to be transitive. (This is why \citeauthor{LadymanLinneboPettigrew:IDPL} \cite[p.\ 23]{LadymanLinneboPettigrew:IDPL} suggest that $\yesrel$ and $\norel$ violate a plausible `minimal requirement' on any notion of indiscernibility.) But it does not immediately follow that one cannot \emph{propose} one of these four grades as a criterion of identity.

\section{Connections to definability theory}\label{s:DefinabilityTheory}
I now want to explore some natural technical questions which have not featured on the radar of philosophers interested in grades of discrimination. These questions concern the relationship between grades of discrimination and \emph{elementary extensions}, and they relate to definability theory. My answers to these questions, together with the Galois Connection of \S\ref{s:Galois}, will yield interesting entailments between the different grades in special cases (to be discussed in \S\ref{s:FinitaryCase}).  To be clear on terminology:
	\begin{define}
		Let $\model{M}$ and $\model{N}$ be $\langmaybe$-structures. Say that $\model{M} \precyes \model{N}$ iff for all $n < \omega$, all $\phi \in \langyes_n$, and all $\overline{e} \in M^n$:
		\begin{align*}
			\model{M} \models \phi(\overline{e}) &\text{ iff }\model{N}\models \phi(\overline{e})  
		\end{align*}
		Say that $\model{M} \precno \model{N}$ iff the above holds with $\langno_n$ in place of $\langyes_{n}$. 
	\end{define}\noindent 
There is a classic result connecting $\langyes$-indiscernibility with the existence of a symmetry in some elementary extension (see e.g.\ \citeauthor{Marker:MT} \cite[Proposition 4.1.5]{Marker:MT}):
\begin{thm}\label{thm:YesSvenStep} 
	For any $\langmaybe$-structure $\model{M}$, the following are equivalent:
	\begin{enumerate}[nolistsep]
		\item $\model{M} \models \phi(\overline{a}) \liff \phi(\overline{b})$, for all $\phi \in \langyes_n$
		\item There is an $\model{N} \succyes \model{M}$ and a symmetry $\pi$ on $\model{N}$ such that $\pi(\overline{a}) = \overline{b}$ 
	\end{enumerate}
\end{thm}\noindent
For present purposes, the immediate import of Theorem \ref{thm:YesSvenStep} is that it yields a new way to characterise $\yesrel$ and $\yesmonad$:
\begin{lem}\label{lem:YesSvenConsequence}
	For any $\langmaybe$-structure $\model{M}$:
	\begin{enumerate}[nolistsep]
		\item\label{yesmonadextend} $a \yesmonad b$ in $\model{M}$ {iff} there is an $\model{N} \succyes \model{M}$ in which $a \symbare b$
 		\item\label{yespairextend} $a \yesrel b$ in $\model{M}$ {iff} there is an $\model{N} \succyes \model{M}$ in which $a \sympair b$ 
	\end{enumerate}
\end{lem}\noindent
This raises a natural question: Is there an $\langno$-analogue of Theorem \ref{thm:YesSvenStep}? There certainly is; but to show this, I need two definitions. First, I need the ordinary notion of a diagram:\footnote{Whilst this notion of diagram invokes `$=$', \citeauthor{Dellunde:EFL} shows that there is a perfectly workable notion of diagram which does not employ `$=$'.}
\begin{define}\label{define:ElemDiagram}
	Let $\langmaybe$ be any signature and $X$ be any set. Then $\langmaybe(X)$ is the signature formed by augmenting $\langmaybe$ with each member of $X$ as a (new) constant. Where $\model{M}$ is an $\langmaybe$-structure, $\edyes(\model{M})$ is the set of $\langyes(M)$-sentences satisfied by the $\langmaybe(M)$-structure formed by letting each $e \in M$ name itself. 
\end{define}\noindent
Next, I need the $\langno$-analogue for a partial elementary map:
\begin{define}
	Let $\model{M}, \model{N}$ be $\langmaybe$-structures. A \emph{proto-correspondence from $\model{M}$ to $\model{N}$} is any relation $\Pi \subseteq M \times N$ such that, for all $n < \omega$, all $\phi \in \langno_n$, and all $\overline{d}\Pi\overline{e}$:
		\begin{align*}
			\model{M} \models \phi(\overline{d})&\text{ iff }\model{N} \models \phi(\overline{e})
		\end{align*}
\end{define}\noindent
So a near-correspondence from $\model{M}$ to $\model{N}$ is a proto-correspondence with domain $M$ and range $N$.  The proof of the $\langno$-analogue of Theorem \ref{thm:YesSvenStep} now amounts to little more than a tweak to \citeauthor{Marker:MT}'s proof of Theorem \ref{thm:YesSvenStep}.\footnote{In more detail: my Lemma \ref{lem:AddOne} tweaks \citeauthor{Marker:MT}'s Lemma 4.16; my Lemma \ref{lem:OneToAll} tweaks \citeauthor{Marker:MT}'s Corollary 4.1.7 (cf.\ also \citeauthor{CasanovasEtAl:EEEFL}'s  \cite{CasanovasEtAl:EEEFL} Lemma 2.7); and my Lemma \ref{lem:SvenoniusCoolDiagram} tweaks \citeauthor{Marker:MT}'s Proposition 4.1.5. The main difference is that I use proto-correspondences rather than partial elementary maps, and in the final step I require a detour, via Lemma \ref{lem:GaloisConsequences}, to obtain a relativity.} I start with two type-realising constructions:
\begin{lem}\label{lem:AddOne}
	Let $\Pi$ be a proto-correspondence from $\model{M}$ to $\model{N}$. For any $a \in M$,  there is as an $\model{O} \succyes \model{N}$ with some $b \in O$ such that $\Pi \cup \{\langle a, b\rangle\}$ is a proto-correspondence from $\model{M}$ to $\model{O}$.
	\begin{proof}
		Define:
		 	\begin{align*}
				\Gamma &= \{\phi(v, \overline{e}) \in \langno_1(\range(\Pi)) \mid \text{for some }n< \omega\text{, some }\phi \in \langno_{n+1}\text{, and  some }\overline{d}\Pi\overline{e},\\
				&\phantom{= \{\phi(v, \overline{e}) \in \langno_1(\range(\Pi)) \mid .}\text{we have }\model{M} \models \phi(a, \overline{d})\}
			\end{align*}
		Consider any $\phi(v, \overline{e}) \in \Gamma$; since $\model{M} \models \exists v \phi(v, \overline{d})$ and $\Pi$ is a proto-correspondence, $\model{N} \models \exists v \phi(v, \overline{e})$. Equally, $\model{N}$ can be treated as a model of $\edyes(\model{N})$. So any finite subset of $\Gamma \cup \edyes(\model{N})$ is satisfiable. Hence, by Compactness, there is a model of $\Gamma \cup \edyes(\model{N})$, which we can regard as $\model{O} \succyes \model{N}$. Now simply let $\Sigma = \Pi \cup \{\langle a, b \rangle\}$, where $\model{O} \models \Gamma(b)$.
	\end{proof}
\end{lem}

\begin{lem}\label{lem:OneToAll}
		Let $\Pi$ be a proto-correspondence from $\model{M}$ to $\model{N}$ with $\model{M} \precyes \model{N}$. Then there is some $\model{O} \succyes \model{N}$ and a proto-correspondence $\Sigma \supseteq \Pi^{-1}$ from $\model{N}$ to $\model{O}$ with $\domain(\Sigma)=N$.
	\begin{proof}
		We construct an elementary chain. Since $\Pi$ is a proto-correspondence and $\model{M} \precyes \model{N}$, we have that for all $\phi \in \langno_n$ and all $\overline{a}\Pi\overline{b}$: 
			\begin{align*}
				\model{N} \models \phi(\overline{b}) \text{ iff }\model{M} \models \phi(\overline{a})\text{ iff }\model{N} \models \phi(\overline{a})
			\end{align*}
		 Defining $\model{O}_0 = \model{N}$ and $\Sigma_0 = \Pi^{-1}$, observe that $\Sigma_0$ is a proto-correspondence from $\model{N}$ to $\model{O}_0$. This is our initial stage in the chain.
			
		Now let $\{e_\alpha \mid \alpha < \kappa\}$ exhaustively enumerate $N$, let $D_\alpha = \domain(\Sigma_0) \cup \{e_\beta  \mid \beta < \alpha\}$ for each $\alpha < \kappa$, and proceed recursively:
	\begin{enumerate}[nolistsep]
		\item[---] Stage $\alpha + 1$: Given a proto-correspondence $\Sigma_\alpha$ from $\model{N}$ to $\model{O}_{\alpha}$ with $\domain(\Sigma_\alpha) = D_\alpha$, use Lemma \ref{lem:AddOne} to obtain an $\model{O}_{\alpha+1} \succyes \model{O}_{\alpha}$ and a proto-correspondence $\Sigma_{\alpha + 1} \supseteq \Sigma_{\alpha}$ from $\model{N}$ to $\model{O}_{\alpha + 1}$ with $\domain(\Sigma_{\alpha+1}) = D_{\alpha+1}$. 
		\item[---] Stage $\alpha$, with $\alpha$ a limit ordinal: let $\model{O}_\alpha = \bigcup_{\beta < \alpha} \model{O}_\beta$ and $\Sigma_\alpha = \bigcup_{\alpha < \beta} \Sigma_\beta$. 
	\end{enumerate}
Now let $\model{O} = \bigcup_{\alpha < \kappa}\model{N}_\alpha$ and $\Sigma = \bigcup_{\alpha < \kappa} \Sigma_\alpha$.
	\end{proof}
\end{lem}
\begin{lem}\label{lem:SvenoniusCoolDiagram}
	Let $\model{M}$ be an $\langmaybe$-structure with $\overline{a}, \overline{b} \in M^n$ such that $\model{M} \models \phi(\overline{a}) \liff \phi(\overline{b})$ for all $\phi \in\langno_n$. Then there is some $\model{N} \succyes \model{M}$ and a near-correspondence $\Pi$ from $\model{N}$ to $\model{N}$ such that $\overline{a}\Pi\overline{b}$.
	\begin{proof}
		Given $\model{M}, \overline{a}, \overline{b}$ as described, we have a proto-correspondence $\Pi_0$ from $\model{M}$ to $\model{M}$ with $\overline{a}\Pi_0\overline{b}$. Setting $\model{M} = \model{M}_0 = \model{N}_0$, we can repeatedly apply Lemma \ref{lem:OneToAll} to construct an elementary chain (solid arrows indicate a proto-correspondence):
	
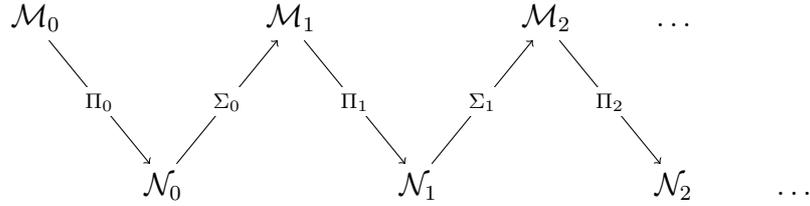
\begin{center}
            \begin{tikzpicture}[description/.style={fill=white,inner sep=2pt}]
            \matrix (m) [matrix of math nodes, row sep=4em,
            column sep=2em, text height=1.5ex, text depth=0.25ex]
            { \model{M}_0 & & \model{M}_1 & & \model{M}_2 & \ldots\\
            & \model{N}_0 & & \model{N}_1 & & \model{N}_2 & \ldots\\ };
            \path[->,font=\scriptsize]
            (m-1-1.300) edge node[description] {$\Pi_0$}(m-2-2.120)
            (m-1-3.300) edge node[description] {$\Pi_1$} (m-2-4.120)
            (m-1-5.300) edge node[description] {$\Pi_2$} (m-2-6.120)
            (m-2-2.60) edge node[description] {$\Sigma_0$} (m-1-3.240)
            (m-2-4.60) edge node[description] {$\Sigma_1$} (m-1-5.240);
		\end{tikzpicture}
	\end{center}
where both $\Sigma_i \subseteq \Pi_i^{-1} \subseteq \Sigma_{i+1}$ and $\model{M}_i \precyes \model{N}_i \precyes \model{M}_{i+1}$ for each $i < \omega$. Define:
	\begin{align*}
		\model{N} &= \bigcup_{i < \omega} \model{N}_i = \bigcup_{i < \omega} \model{M}_i\\
		\Pi &= \bigcup_{i < \omega} \Pi_i
	\end{align*}
It is routine to check that $\Pi$ and $\model{N}$ have the required properties. 
\end{proof}
\end{lem}\noindent
We can now obtain our $\langno$-analogue of Theorem \ref{thm:YesSvenStep}:
\begin{thm}\label{thm:NoSvenStep}For any $\langmaybe$-structure $\model{M}$, the following are equivalent:
	\begin{enumerate}[nolistsep]
		\item\label{Svenonius:Map} $\model{M} \models \phi(\overline{a}) \liff \phi(\overline{b})$, for all every $\phi \in \langno_n$
		\item\label{Svenonius:RelativityYes} There is an $\model{N} \succyes \model{M}$ and a relativity $\Pi$ on $\model{N}$ such that $\overline{a}\Pi\overline{b}$
		\item\label{Svenonius:RelativityNo} There is an $\model{N} \succno \model{M}$ and a relativity $\Pi$ on $\model{N}$ such that $\overline{a}\Pi\overline{b}$ 
	\end{enumerate}
		\begin{proof}(\ref{Svenonius:Map}) $\Rightarrow$ (\ref{Svenonius:RelativityYes}). Use Lemma \ref{lem:SvenoniusCoolDiagram} to obtain a near-correspondence, then use Lemma \ref{lem:GaloisConsequences} to extend this to a relativity.
		
		(\ref{Svenonius:RelativityYes}) $\Rightarrow$ (\ref{Svenonius:RelativityNo}). Trivial. 
		
		(\ref{Svenonius:RelativityNo}) $\Rightarrow$ (\ref{Svenonius:Map}). $\Pi$ is a relativity, and hence a near-correspondence by Lemma \ref{lem:RelativenessPreservation}; the result now follows since $\model{M} \precno\model{N}$.	
	\end{proof}
\end{thm}\noindent
This Theorem lends yet more weight to the claim that relativeness correspondences are the $\langno$-analogue of symmetries. Moreover, it immediately yields a new way to characterise $\norel$ and $\nomonad$ (compare Lemma \ref{lem:YesSvenConsequence}):
\begin{lem}\label{lem:NoSvenConsequence} For any $\langmaybe$-structure $\model{M}$:
	\begin{enumerate}[nolistsep]
		\item $a \nomonad b$ in $\model{M}$ \text{iff} there is an $\model{N} \succyes \model{M}$ in which $a \relativebare b$ 
		\item $a \norel b$ in $\model{M}$ \text{iff} there is an $\model{N} \succyes \model{M}$ in which $a \relativepair b$
	\end{enumerate}\noindent
	Moreover, both claims hold with $\succno$ in place of $\succyes$. 		 
\end{lem}\noindent
Before continuing with the main aims of this paper, it is worth briefly stopping to smell the roses. Theorem \ref{thm:YesSvenStep} is sometimes used as a stepping stone to the following foundational result of definability theory (notation clarified in endnote):\footnote{\citeauthor{Beth:PMTD} \cite{Beth:PMTD} proved (\ref{YesSven:definable}) $\Leftrightarrow$ (\ref{YesSven:Beth}); \citeauthor{Svenonius:TPM} \cite{Svenonius:TPM} proved (\ref{YesSven:definable}) $\Leftrightarrow$ (\ref{YesSven:invariance}). $(\model{M}, U)$ is the $\langmaybe\mathord{\cup}\{R\}$-structure formed from $\model{M}$ by allowing $R$ to pick out $U$. As one would expect, $\pi(V) = \{\pi(\overline{e}) \mid \overline{e} \in V\}$, and equally $\Pi(V) = \{\overline{e} \mid \text{there are }\overline{d} \in V\text{ such that }\overline{d}\Pi\overline{e}\}$.}
\begin{thm}[Beth--Svenonius Theorem, $\langyes$-case]\label{thm:YesSven}
For any $\langmaybe$-structure $\model{M}$ with $R \notin \langmaybe$ and $U \subseteq M^n$, the following are equivalent:
	\begin{enumerate}[nolistsep]
		\item\label{YesSven:definable} $(\model{M}, U) \models \forall \overline{v}(\phi(\overline{v}) \liff R\overline{v})$ for some $\phi \in \langyes_n$.
		\item\label{YesSven:invariance} For every $(\model{N}, V) \succyes (\model{M}, U)$ and every symmetry $\pi$ of $\model{N}$: $V = \pi(V)$. 
		\item\label{YesSven:Beth} For any $\langmaybe$-structure $\model{N}$ and any sets $V_0, V_1 \subseteq N^n$: if $(\model{N}, V_0)$, $(\model{N}, V_1)$ and $(\model{M}, U)$ all satisfy the same $\langyes$-sentences, then $V_0 = V_1$. 
	\end{enumerate}
\end{thm}\noindent
Pleasingly, we can use Theorem \ref{thm:NoSvenStep} as a stepping-stone to an $\langno$-analogue of this result; indeed, the main steps are exactly as in the $\langyes$-case:\footnote{The only difficult step in either Theorem is (\ref{NoSven:invarianceYes}) $\Rightarrow$ (\ref{NoSven:definable}). For the $\langyes$-case, see e.g.\ \citeauthor{Poizat:MT} \cite[Proposition 9.2]{Poizat:MT}. To prove the $\langno$-case, we simply tweak Poizat's proof by invoking Theorem \ref{thm:NoSvenStep} rather than Theorem \ref{thm:YesSvenStep}, and considering $n$-types in the sense of $\langno_n$, rather than $\langyes_n$. \citeauthor{Dellunde:EFSM} \cite[p.\ 5]{Dellunde:EFSM} and \citeauthor{KeislerMiller:CWE} \cite[p.\ 3]{KeislerMiller:CWE} have shown that $\langno_n$-types behave as one would hope.}
\begin{thm}[Beth--Svenonius Theorem, $\langyes$-case]\label{thm:NoSven}
For any $\langmaybe$-structure $\model{M}$ with $R \notin \langmaybe$ and $U \subseteq M^n$, the following are equivalent:
	\begin{enumerate}[nolistsep]
		\item\label{NoSven:definable} $(\model{M}, U) \models \forall \overline{v}(\phi(\overline{v}) \liff R\overline{v})$ for some $\phi \in \langno_n$. 
		\item\label{NoSven:invarianceYes} For every $(\model{N}, V) \succyes (\model{M}, U)$ and every relativity $\Pi$ of $\model{N}$: $V = \Pi(V)$
		\item\label{NoSven:invarianceNo} For every $(\model{N}, V) \succno (\model{M}, U)$ and every relativity $\Pi$ of $\model{N}$: $V = \Pi(V)$
		\item\label{NoSven:Beth} For any $\langmaybe$-structure $\model{N}$ and any sets $V_0, V_1 \subseteq N^n$: if $(\model{N}, V_0)$, $(\model{N}, V_1)$ and $(\model{M}, U)$ all satisfy the same $\langno$-sentences, then $V_0 = V_1$. 
	\end{enumerate}
\end{thm}\noindent

\section{Entailments in the finitary case}\label{s:FinitaryCase}
The results of the previous section immediately yield a special case of Theorem \ref{thm:BigMap}, obtained by restricting our attention to \emph{finitary} structures.\footnote{An alternative route to Theorem \ref{thm:FiniteMap} merits comment. We can use finitary isomorphism systems to prove the coincidence of grades of $\langyes$-indiscernibility with grades of symmetry in finite structures, without invoking the results from \S\ref{s:DefinabilityTheory}. (For an introduction to finitary isomorphism systems, see \citeauthor{EbbinghausEtAl:ML} \cite[chapter XI]{EbbinghausEtAl:ML}.) \citeauthor{CasanovasEtAl:EEEFL} \cite[Definitions 4.1--4.2]{CasanovasEtAl:EEEFL} define the $\langno$-analogue of finitary isomorphism systems. It turns out that we can use these, analogously, to prove the coincidence of grades of $\langno$-indiscernibility with grades of relativity in finite structures.} This special case has already attracted some attention, since it is philosophically interesting,\footnote{\citeauthor{CaultonButterfield:KILM} \cite[Theorem 2]{CaultonButterfield:KILM} prove a special case of the mutual entailment between $\sympair$ and $\yesrel$, and $\symbare$ and $\yesmonad$, on the assumption that $\langmaybe$ is finite and relational. \citeauthor{Ketland:II}'s \cite*[p.\ 2]{Ketland:II} attention is entirely restricted to finite relational signatures. \citeauthor{LinneboMuller:WD} \cite[Theorem 3]{LinneboMuller:WD} note that \emph{witness-discernibility} (a further notion, which I have not discussed) is equivalent to $\norel$ in finite structures, and outline several reasons for focussing on finite structures.} and the following result completes the picture.
\begin{thm}[Entailments between the grades, finite structures]\label{thm:FiniteMap}
These Hasse Diagrams characterise the entailments between our grades of discrimination, when we restrict our attention to finite structures:
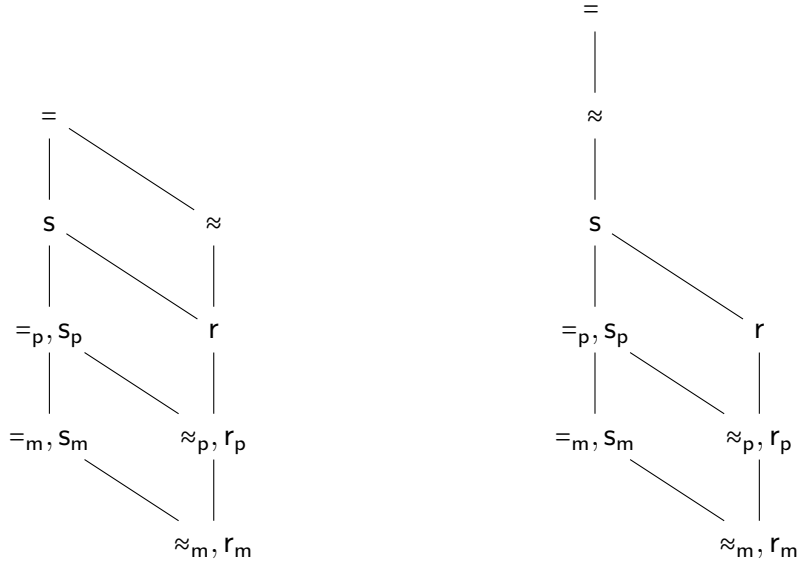
\begin{quote}
	\begin{tikzpicture}
		\matrix (m) [matrix of math nodes, row sep=2em, 
		column sep=2em, text height=1.5ex, text depth=0.25ex] 
		{  &  & \\
		  = & \\
		 \symtotal & \foi \\
		\yesrel, \sympair & 		 \relativetotal \\
		 \yesmonad, \symbare & \norel, \relativepair \\
		 & \nomonad, \relativebare \\};
		\draw(m-2-1)--(m-3-1)--(m-4-1)--(m-5-1)--(m-6-2);
		\draw(m-2-1)--(m-3-2)--(m-4-2)--(m-5-2)--(m-6-2);
		\draw(m-3-1)--(m-4-2);
		\draw(m-4-1)--(m-5-2);
		\end{tikzpicture}
\hfill
	\begin{tikzpicture}
		\matrix (m) [matrix of math nodes, row sep=2em, 
		column sep=2em, text height=1.5ex, text depth=0.25ex] 
		{ = &  & \\
		  \foi & \\
		 \symtotal &\\
		\yesrel, \sympair & 		 \relativetotal \\
		 \yesmonad, \symbare & \norel, \relativepair \\
		 & \nomonad, \relativebare \\};
		\draw(m-1-1)--(m-2-1)--(m-3-1)--(m-4-1)--(m-5-1)--(m-6-2);
		\draw(m-4-2)--(m-5-2)--(m-6-2);
		\draw(m-3-1)--(m-4-2);
		\draw(m-4-1)--(m-5-2);
		\end{tikzpicture}
\end{quote}
	The left diagram is restricted to finite structures with arbitrary signatures; the right diagram is restricted to finite structures with \emph{relational} signatures.
	\begin{proof}
		Most of this is supplied by Theorem \ref{thm:BigMap}. For the remainder observe that if $\model{M}$ is finite, then $\model{N} \succyes \model{M}$ iff $\model{M} = \model{N}$. It follows from Lemma \ref{lem:YesSvenConsequence} that $\yesrel$ entails $\sympair$ and that $\yesmonad$ entails $\symbare$; and similarly with Lemma \ref{lem:NoSvenConsequence}. 
	\end{proof}
\end{thm}\noindent
However, a little work will yield an even stronger result: the grades of relativity and the grades of $\langno$-indiscernibility also entail each other when the structure's $\foi$-quotient is finite. To show this, we need a few results. The first tells us when $\foi$ is definable in a structure:\footnote{For the case where $\langmaybe$ is finite and relational, see \citeauthor{Ketland:II} \cite[Definition 2.3]{Ketland:II}.}
	\begin{lem}\label{lem:QuotientBehaviour}Let $\model{M}$ be any $\langmaybe$-structure. If either $\langmaybe$ is finite and relational or $\quotient{M}$ is finite, then there is an $\langno_2$-formula which defines $\foi$ in $\model{M}$. However, the restrictions are necessary.
		\begin{proof}
			\emph{Case when $\langmaybe$ is finite and relational.} Let $\phi_1, \ldots, \phi_n$ enumerate all the atomic $\langmaybe$-formulas. By Lemma \ref{lem:AlternateCharacterisationOfFOI}, the following $\langno_2$-formula defines $\foi$ in $\model{M}$:
		\begin{align*}
			\bigwedge_{i = 1}^n \forall \overline{v} (\phi_i (x,\overline{v}) \liff \phi_i(y, \overline{v}))
		\end{align*}		
		
			\emph{Case when $\quotient{M}$ is finite.} Let $\quotient{e_1}, \ldots, \quotient{e_m}$ exhaustively enumerate $\quotient{M}$ without repetition. So for all $i \neq j$ between $1$ and $m$, we have ${e_i} \foi {e_j}$; hence by Lemma \ref{lem:AlternateCharacterisationOfFOI} there is some $\phi_{i,j} \in \langno_2$ such that $\model{M} \nmodels \phi_{i,j}({e_i}, {e_i}) \liff \phi_{i,j}({e_i}, {e_j})$, and so the following $\langno_2$-formula defines $\foi$ in $\model{M}$:
			\begin{align*}
				\bigwedge_{i \neq j}(\phi_{i,j}(x, x) \liff \phi_{i,j}(x, y))
			\end{align*}
			
		\emph{The necessity of the restrictions.} Let $\langmaybe$ contain one-place predicates $P_i$ for all $i < \omega$, and a single two-place predicate $R$. Define:
	\begin{align*}
		\NotDefineIdentity :=& \mathbb{N}\\
		P_i^\model{\NotDefineIdentity} := & \{6i, 6i+1\}, \text{for all }0 < i < \omega\\
		R^\model{\NotDefineIdentity} := & 
		\{\langle 2, 1\rangle, \langle 1, 0\rangle\} \cup \{\langle 2, 6n -2\rangle, \langle 6n-2, 6n\rangle, \langle 6n-2, 6n+2\rangle \mid 0 < n < \omega\}\ \cup \\
	&	\{\langle 3, 6n -1\rangle, \langle 6n-1, 6n+1\rangle, \langle 6n-1, 6n+3\rangle \mid 0 < n < \omega\}
	\end{align*}
We can represent $\model{\NotDefineIdentity}$ more perspicuously as follows:
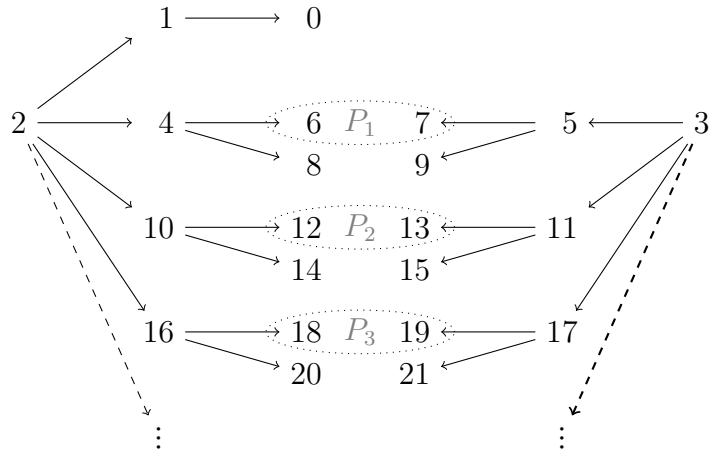
\begin{quote}\centering
\begin{tikzpicture}
		\node(a) {$2$};
		\node[right=3em of a](a1) {$\phantom{1}4$};
		\node[right=3em of a1](a11) {$\phantom{1}6$};
		\node[below=0em of a11](a10) {$\phantom{1}8$};
		\path[->] (a1) edge (a11) edge (a10);
		\node[right=0em of a11](P1) {\textcolor{gray}{$P_1$}};
		\node[below=2em of a1](a2) {$10$};
		\node[right=3em of a2](a21) {$12$};
		\node[below=0em of a21](a20) {$14$};
		\path[->] (a2) edge (a21) edge (a20);
		\node[right=0em of a21](P2) {\textcolor{gray}{$P_2$}};
		\node[below=2em of a2](a3) {$16$};
		\node[right=3em of a3](a31) {$18$};
		\node[below=0em of a31](a30) {$20$};
		\path[->] (a3) edge (a31) edge (a30);
		\node[right=0em of a31](P3) {\textcolor{gray}{$P_3$}};
		\node[below=2em of a3](adots) {$\vdots$};
		\node[above=2em of a1](a0) {$\phantom{1}1$};
		\node[right=3em of a0](a00) {$\phantom{1}0$};
		\draw[dashed, ->] (a)--(adots);
		\draw[->] (a0)--(a00);
		\path[->] (a) edge (a1) edge (a2) edge (a3) edge (a0);		
		\node[right=0em of P1](b11) {$\phantom{2}7$};
		\node[right=3em of b11](b1) {$\phantom{1}5$};
		\node[right=3em of b1](b) {$3$};
		\node[below=0em of b11](b10) {$\phantom{2}9$};
		\draw[->] (b)--(b1);
		\path[->] (b1) edge (b11) edge (b10);
		\node[right=0em of P2](b21) {$13$};
		\node[right=3em of b21](b2) {$11$};
		\node[below=0em of b21](b20) {$15$};
		\node[right=0em of P3](b31) {$19$};
		\node[right=3em of b31](b3) {$17$};
		\node[below=0em of b31](b30) {$21$};
		\node[below=2em of b3](bdots) {$\vdots$};
		\draw[->] (b)--(b2);
		\path[->] (b2) edge (b21) edge (b20);
		\draw[->] (b)--(b3);
		\path[->] (b3) edge (b31) edge (b30);
		\draw[->][dashed, thick] (b)--(bdots);
		\draw[dotted] (P1) ellipse (3em and 0.7em);
		\draw[dotted] (P2) ellipse (3em and 0.7em);
		\draw[dotted] (P3) ellipse (3em and 0.7em);
	\end{tikzpicture}
\end{quote}
	I claim that $\quotient{\model{\NotDefineIdentity}}\models \phi(\quotient{2}) \liff \phi(\quotient{3})$ for all $\phi \in \langno_1$. To prove this, fix $\phi \in \langno_1$ and let $\lang{K}$ be the (necessarily finite) set of $\langmaybe$-predicates appearing in $\phi$. Where $\model{\NotDefineIdentity}^*$ is the $\lang{K}$-reduct of $\model{\NotDefineIdentity}$, we have $\quotient{2} \sympair \quotient{3}$ in $\quotient{\model{\NotDefineIdentity}^*}$, and hence $2 \relativepair 3$ in $\model{\NotDefineIdentity}^*$ by Lemma \ref{lem:QuotientSymHom}. Lemma \ref {lem:RelativenessPreservation} now yields that $\model{H}^* \models \phi(2) \liff \phi(3)$, and hence $\model{H} \models \phi(2) \liff \phi(3)$. Now apply Lemma \ref{lem:Quotient}.
	
	However, where $\psi \in \langyes_1$ abbreviates:
	\begin{align*}
		\exists x (Rvx \land \forall y \forall z ((Rxy \land Rxz) \lonlyif y = z)))
	\end{align*}
	we have  $\quotient{\model{\NotDefineIdentity}} \models \psi(\quotient{2}) \land \lnot \psi(\quotient{3})$. So no $\langno_2$-formula can define $=$ in $\quotient{\model{\NotDefineIdentity}}$; and hence no $\langno_2$-formula can define $\foi$ in $\model{\NotDefineIdentity}$, by Lemma \ref{lem:Quotient}.
	\end{proof}
\end{lem}\noindent
We already knew that $a \foi b$ in $\model{M}$ iff $\quotient{a} = \quotient{b}$ in $\quotient{\model{M}}$. Lemma \ref{lem:QuotientBehaviour} allows us, under special circumstances, to obtain analogous results for our other grades of indiscernibility:
\begin{lem}\label{lem:IndiscernibilityReducts}
Let $\model{M}$ be any $\langmaybe$-structure. If either $\langmaybe$ is finite and relational or $\quotient{M}$ is finite, then:
		\begin{enumerate}[nolistsep]
		\item $a \norel b$ in $\model{M}$ iff $\quotient{a} \yesrel \quotient{b}$ in $\quotient{\model{M}}$
		\item $a \nomonad b$ in $\model{M}$ iff $\quotient{a} \yesmonad \quotient{b}$ in $\quotient{\model{M}}$ 
	\end{enumerate}
	\begin{proof}
		By Lemma \ref{lem:QuotientBehaviour}, given either assumption, some $\langno_2$-formula defines $\foi$ in $\model{M}$. The same formula defines $=$ in $\quotient{\model{M}}$, and the result follows via Lemma \ref{lem:Quotient}.
	\end{proof}
\end{lem}\noindent
Finally, the Galois Connection of \S\ref{s:Galois} allows us to extend Theorem \ref{thm:FiniteMap}, as desired:
\begin{lem}\label{lem:FiniteMapRelativeHalf}
	For structures with finite $\foi$-quotients:
	\begin{enumerate}[nolistsep]
		\item\label{Finite:RelativePair} $\norel$ entails $\relativepair$, and vice versa
		\item\label{Finite:RelativeBare} $\nomonad$ entails $\relativebare$, and vice versa.
	\end{enumerate}
	\begin{proof}
	Combine Theorem \ref{thm:FiniteMap} with Lemmas \ref{lem:IndiscernibilityReducts} and \ref{lem:QuotientSymHom}.
		\end{proof}
\end{lem}\noindent
Note that Theorem \ref{thm:FiniteMap} does not hold for grades of $\langyes$-discernibility/symmetry. To see this, let $\model{\TwoBallBig}^*$ be a superstructure of $\model{\TwoBallBig}$ obtained by making $R$ reflexive. Whilst $\quotient{\model{\TwoBallBig}^*}$ has only two members and its signature is finite and relational, no symmetry on $\model{\TwoBallBig}^*$ sends $1$ to any element in $\mathbb{R} \setminus \mathbb{N}$.

\section{Capturing grades of discrimination}\label{s:Capturing}
All twelve grades of discrimination have fairly straightforward definitions. However, the grades of indiscernibility are defined in terms of satisfaction of object-language formulas, whereas the grades of symmetry and relativity are defined are defined metalinguistically. It is natural to ask whether this is essential. More precisely, I shall ask which of the grades are \emph{capturable}, in the following sense:
\begin{define}
	For any $\langmaybe$-structure $\model{M}$, say that $\Gamma \subseteq \langmaybe_2$ \emph{captures} $\mathrm{R}$ \emph{in $\model{M}$} iff:
		\begin{align*}
			\text{for all }a, b \in M: a\mathrm{R}b\text{ in }\model{M}&\text{ iff }\model{M} \models \phi(a,b)\text{ for every }\phi \in \Gamma
		\end{align*}
	Say that $\mathrm{R}$ is \emph{capturable$\yesstroke$ in $\model{M}$} iff some $\Gamma \subseteq \langyes_{2}$ captures $\mathrm{R}$ in $\model{M}$ (similarly for \emph{capturable$\nostroke$}). Say that $\mathrm{R}$ is \emph{universally capturable$\yesstroke$} \text{iff} some single $\Gamma \subseteq \langyes_{2}$ captures $\mathrm{R}$ in every $\langmaybe$-structure (similarly for \emph{universally capturable$\nostroke$}).

\end{define}\noindent 
I shall consider capturability for each of the three families of grades of discrimination, starting with the grades of indiscernibility:
\begin{lem}\label{lem:CaptureIndiscernible}
		\begin{enumerate}[nolistsep]
			\item\label{CaptureIndiscernible:yes} $=$, $\yesrel$, $\yesmonad$ are universally capturable$\yesstroke$
			\item\label{CaptureIndiscernible:no} $\foi$, $\norel$, $\nomonad$ are universally capturable$\nostroke$
			\item\label{CaptureIndiscernible:nonot} There is a structure in which none of $=$, $\yesrel$ and $\yesmonad$ is capturable$\nostroke$
		\end{enumerate}
	\begin{proof}
		(\ref{CaptureIndiscernible:yes}) and (\ref{CaptureIndiscernible:no}). Obvious. 
				
		(\ref{CaptureIndiscernible:nonot}). Let $\model{\ThreeThreeSphere}$ be the following graph:
			\begin{quote}
			\centering
			\begin{tikzpicture}[descr/.style={fill=white,inner sep=2.5pt}] 
			\matrix (m) [matrix of math nodes, row sep=0em, 
			column sep=2em, text height=1.5ex, text depth=0.25ex] 
			{\phantom{c} & 1 & 2 &3  \\
			& 4 & 5 & 6\\ 
			& 7 & 8 & 9 \\ };
			\draw[-, thick] (m-1-2)--(m-1-3);
			\draw[-, thick] (m-2-2)--(m-2-3);
			\draw[-, thick] (m-3-2)--(m-3-3)--(m-3-4);
			\end{tikzpicture} 
			\end{quote}
		For all $\phi \in \langno_{2}$, $\model{\ThreeThreeSphere} \models \phi(1, 4) \liff \phi(1, 7)$ even though $1 \yesrel 4$ and $1 \nyesmonad 7$.
	\end{proof}
\end{lem}
\begin{lem}\label{lem:CaptureSymmetries}
	\begin{enumerate}[nolistsep]
		\item\label{CaptureSymTotal} $\symtotal$ is universally capturable$\yesstroke$
		\item\label{spsmyesfinite} For finite structures: $\sympair$ and $\symbare$ are universally capturable$\yesstroke$
		\item \label{spsmyes} There is a structure in which $\sympair$ and $\symbare$ are not capturable$\yesstroke$
		\item\label{SymmetryUncapturable}
There is a finite structure in which none of $\symtotal$, $\sympair$ and $\symbare$ is capturable$\nostroke$
	\end{enumerate}
	\begin{proof}
		(\ref{CaptureSymTotal}). For each atomic $\phi \in \langyes_{n+2}$, define:
		\begin{align*}
			x \simeq_\phi y &:= \forall \overline{v} \left(\bigwedge_{i=1}^n(v_i  \neq x \land v_i \neq y) \lonlyif \left[\phi (x, y, \overline{v}) \liff \phi(y, x, \overline{v}) \right]\right)
		\end{align*}
		By Lemma \ref{lem:IsomorphismPreservation}, $\symtotal$ is universally captured$\yesstroke$ by the set of all such $\simeq_\phi$.
	
		(\ref{spsmyesfinite}). From Theorem \ref{thm:FiniteMap} and Lemma \ref{lem:CaptureIndiscernible}.

		(\ref{spsmyes}). Let $\model{\ThreeBallBig}$ comprise two disjoint copies of the complete countable graph, with a disjoint copy of a complete uncountable graph, i.e.:
			\begin{align*}
				\ThreeBallBig := & \mathbb{R}\\
				R^{\model{\ThreeBallBig}} :=& \{\langle m, n\rangle \in \mathbb{N}^2\mid m \neq n\text{ and }m+n\text{ is even}\}  \cup \{\langle p, q \rangle \in (\mathbb{R} \setminus \mathbb{N})^2 \mid p \neq q\}
			\end{align*}
		By taking a Skolem Hull containing $1, 2$ and some $e \in \mathbb{R} \setminus \mathbb{N}$, it is clear that:
			\begin{align*}
				\model{\ThreeBallBig} \models \phi(1, 2) \liff \phi(1,e)
			\end{align*}
		for any $\phi \in \langyes_2$. However, $1 \sympair 2$ in $\model{\ThreeBallBig}$, whereas $1 \nsymbare e$ in $\model{\ThreeBallBig}$. 

		(\ref{SymmetryUncapturable}) In $\model{\ThreeThreeSphere}$ from Lemma \ref{lem:CaptureIndiscernible}, $1 \symtotal 2$, whereas $7 \nsymbare 8$. However, for all $\phi \in \langno_{2}$, $\model{\quotient{\ThreeThreeSphere}} \models \phi(\quotient{1}, \quotient{2}) \liff \phi(\quotient{7}, \quotient{8})$ and hence $\model{\ThreeThreeSphere} \models \phi(1, 2) \liff \phi(7, 8)$.
	\end{proof}
\end{lem}
\begin{lem}\label{lem:RelativityCapturable}
	\begin{enumerate}[nolistsep]
		\item\label{relativetotalcapturable} $\relativetotal$ is universally capturable$\nostroke$
		\item\label{rprmno} For structures with finite $\foi$-quotients: $\relativepair$ and $\relativebare$ are universally cap\-tu\-ra\-ble$\nostroke$
		\item\label{rprmyes} There is a structure in which neither of $\relativepair$ and $\relativebare$ is capturable$\yesstroke$
	\end{enumerate}
	\begin{proof}
		(\ref{relativetotalcapturable}). Let $\Gamma$ be the set of all $\langno_2$-formulas of the form:
			\begin{align*}
				\forall \overline{v} \left(\bigwedge_{i=1}^{n}\left[\phi_{i}(x,x) \land \lnot \phi_{i}(x, v_{i}) \land \psi_{i}(y,y) \land \lnot \psi_{i}(y, v_{i})\right] \lonlyif \left[\theta(x, y, \overline{v}) \liff \theta(y, x, \overline{v})\right]\right)
			\end{align*}
		for any $n < \omega$, any $\phi_1, \ldots, \phi_n, \psi_1, \ldots, \psi_n \in \langno_2$, and any $\theta \in \langno_{n+2}$. I claim that $\Gamma$ captures $\relativetotal$ in any $\langmaybe$-structure $\model{M}$.
		
		First, suppose $a \relativetotal b$ in $\model{M}$. Fix some $\gamma \in \Gamma$, and some $\overline{e} \in M^n$. Suppose that:
		$$\model{M} \models \bigwedge^n_{i=1}\left[\phi_{i}(a,a) \land \lnot \phi_{i}(a, e_{i}) \land \psi_{i}(b,b) \land \lnot \psi_{i}(b, e_{i})\right]$$
		Then by Lemma \ref{lem:AlternateCharacterisationOfFOI}, $e_{i} \nfoi a$ and $e_{i} \nfoi b$ for each $1 \leq i \leq n$. Since $a \relativetotal b$, Lemma \ref{lem:RelativenessPreservation} tells us that $\model{M} \models \theta(a,b,\overline{e}) \liff \theta(b, a, \overline{e})$. Hence $\model{M} \models \gamma(a, b)$, for any $\gamma \in \Gamma$.
		
		 Next, suppose $\model{M} \models \gamma(a,b)$, for all $\gamma \in \Gamma$. I claim that the following is a near-correspondence from $\model{M}$ to $\model{M}$:
		 $$\Pi = \{\langle a, b\rangle, \langle b, a\rangle\} \cup \{\langle x, x\rangle \mid x \nfoi a \text{ and }x\nfoi b\}$$
		To show this, fix $n < \omega$, $\theta \in \langno_{n+2}$ and $\overline{e} \in M^n$ such that $e_i \nfoi a$ and $e_{i} \nfoi b$ for each $1 \leq i \leq n$. Since each $e_{i} \nfoi a$ and $e_{i} \nfoi b$, by Lemma \ref{lem:AlternateCharacterisationOfFOI} there are formulas $\phi_i, \psi_i \in \langno_2$ for each $1 \leq i \leq n$ such that ${M} \models \phi_{i}(a, a) \land \lnot\phi_{i}(a, e_{i})$ and $\model{M} \models \psi_{i}(b, b) \land \lnot \psi_{i}(b, e_{i})$. Conjoining these, we get:
		\begin{align*}
				\model{M} \models & \bigwedge_{i=1}^{n}\left[\phi_{i}(a,a) \land \lnot \phi_{i}(a, e_{i}) \land \psi_{i}(b,b) \land \lnot \psi_{i}(b, e_{i})\right]
			\end{align*}
		Since $\model{M} \models \gamma(a, b)$ for all $\gamma \in \Gamma$, we obtain that, for all $\theta \in \langno_{n+2}$: 
		$$\model{M} \models \theta(a, b, \overline{e}) \liff \theta(b, a, \overline{e})$$
		Generalising, $\Pi$ is a near-correspondence. By the Galois Connection of Theorem \ref{thm:GaloisConnection}, $(\Pi^{\galdown})^{\galup}$ is a relativity on $\model{M}$; and so $a \relativetotal b$.
			
		(\ref{rprmno}). From Lemmas \ref{lem:FiniteMapRelativeHalf} and  \ref{lem:CaptureIndiscernible}.

		(\ref{rprmyes}). Exactly as in Lemma \ref{lem:CaptureSymmetries}, case (\ref{spsmyes}).
		\end{proof}
	\end{lem}\noindent
Lemmas \ref{lem:CaptureIndiscernible}--\ref{lem:RelativityCapturable} can be summarised as follows:
\begin{thm}[Capturing the grades]\label{thm:Capturability} The following table exhaustively details the capturability of each grade of discrimination:
\begin{table}[h]
\setlength\tabcolsep{0.025\textwidth}
\begin{tabular}{@{}p{0.2\textwidth}p{0.2\textwidth}p{0.2\textwidth}@{}}
			Grade & Capturable$\yesstroke$ & Capturable$\nostroke$\\
			\hline
			$=$ & \verified & \falsified \\
			$\yesrel$ & \verified & \falsified \\
			$\yesmonad$ & \verified & \falsified  \\\\
			$\foi$ & \verified & \verified \\
			$\norel$ & \verified & \verified \\
			$\nomonad$ & \verified & \verified \\\\
			$\symtotal$& \verified & \falsified \\
			$\sympair$& \finitecase & \falsified \\
			$\symbare$& \finitecase & \falsified \\\\
			$\relativetotal$& \verified & \verified \\
			$\relativepair$& \finitequotient & \finitequotient \\
			$\relativebare$& \finitequotient & \finitequotient 
	\end{tabular}
\end{table}
\end{thm}\noindent 
The table of Theorem \ref{thm:Capturability} should be read with the following key:
		\begin{enumerate}[nolistsep]
		\item[\verified] universally capturable
		\item[\falsified] there is an $\langmaybe$-structure in which the grade is not capturable
		\item[\finitecase\ ] universally capturable when we restrict attention to finite structures; but there are counterexamples elsewhere
		\item[\finitequotient\ ] universally capturable when we restrict attention to structures with finite $\foi$-quotients; but there are counterexamples elsewhere  
	\end{enumerate} 
This demonstrates, once again,  that grades of $\langyes$-discernibility are to grades of symmetry, as grades of $\langno$-discernibility are to grades of relativity. More interestingly, though, Theorem \ref{thm:Capturability} bears directly upon the philosophical search for reductive criteria of identity. 

As mentioned in \S\ref{s:Preliminaries}, much of the interest in grades of discrimination comes from their potential to provide us with a criterion of identity, possibly a reductive one. However, if a grade of discrimination cannot be captured by some set of formulas in the object language, this should bar it from use in any reductive criterion of identity. After all, if the grade must be invoked as a \emph{primitive} at the level of the object language, it is unclear why we should not simply allow ourselves to take {identity} itself as a primitive in the object language. The situation will be no better, in this regard, if the grade can only be captured$\yesstroke$ and not captured$\nostroke$. Consequently, no grade of $\langyes$-indiscernibility or symmetry can provide a \emph{reductive} criterion of identity. 

The remaining candidates for reductive criteria of identity are therefore the grades of $\langno$-indiscernibility and relativity. However, in the special cases when they are capturable$\nostroke$---which we require if we seek a reductive criterion of identity---two of the grades of $\langno$-indiscernibility are simply co-extensive with two of the grades of $\langno$-indiscernibility (see Theorem \ref{thm:FiniteMap}). Hence the only plausible distinct candidates for a reductive criterion of identity are, in order of entailment: $\foi$, $\relativetotal$, $\norel$ and $\nomonad$.

This does not show, though, that the remaining grades of discrimination are philosophically uninteresting. After all, we might simply be interested in providing an illuminating but \emph{non-reductive} answer to the general question: \emph{When are objects identical?} To repeat an example from \S\ref{s:Preliminaries}: if we have become convinced that nature abhors a (non-trivial) symmetry, then $\symbare$ could serve as a non-reductive, non-trivial criterion of identity, even though it is uncapturable$\yesstroke$. 

\section{Symmetry in all elementary extensions}\label{s:AllElementaryExtensions}
In \S\ref{s:DefinabilityTheory}, I connected the grades of indiscernibility with the existence of a symmetry/relativity in \emph{some} elementary extensions. To close this paper, I wish to consider what happens when we consider the existence of a symmetry or relativity in \emph{all} elementary extensions. In particular, I shall demonstrate a neat connection between $\foi$ and symmetries in elementary extensions. To show this, I first require a general method for constructing such elementary extensions:\footnote{\citeauthor{Monk:ML} \cite[Theorem 29.16]{Monk:ML} described this explicitly; \citeauthor{Grzegorczyk:CC} \cite[p.\ 41]{Grzegorczyk:CC} earlier mentioned it in passing, implying it was mathematical folklore. The method was rediscovered by philosophers, e.g.\ \citeauthor{Ketland:II} \cite[p.\ 7]{Ketland:II} and \citeauthor{LadymanLinneboPettigrew:IDPL} \cite[Theorem 8.14]{LadymanLinneboPettigrew:IDPL}. However, all these authors restrict their attentions to relational signatures.}
\begin{lem}\label{lem:Inflate}
	Let $\model{M}$ be an $\langmaybe$-structure with $a \in M$, and let $D$ be a set such that $M\cap D = \emptyset$. Then there is an $\langmaybe$-structure $\model{N} \succno \model{M}$ with $N = M \cup D$, such that $a \foi d$ in $\model{N}$ for all $d \in D$.
	\begin{proof}
		Define $\sigma : N \longrightarrow M$ by: $\sigma(x) = x$ if $x \in M$, and $\sigma(d) = a$ if $d \in D$. Set:
		\begin{align*}
			R^\model{N} &= \{\overline{e} \in N^n \mid \sigma(\overline{e}) \in R^\model{M} \} & & \text{all }n\text{-place }\langmaybe\text{-predicates }R\\
			f^\model{N}(\overline{e}) &= f^\model{M}(\sigma(\overline{e})) & & \text{all }n\text{-place }\langmaybe\text{-function-symbols }f\text{ and all }\overline{e} \in N^n
		\end{align*}
	I claim that, for each $\langmaybe$-term $\tau$, all $\overline{d} \in D^{m}$ and all $\overline{e} \in M^{n}$:
		\begin{align*}
			\tau^{\model{N}}(\overline{d}, \overline{e}) &= \tau^{\model{N}}(\overline{a}, \overline{e}) = \tau^{\model{M}}(\overline{a}, \overline{e})
		\end{align*}
	(where $a_i = a$ for all $1 \leq i \leq m$). This is proved by induction on complexity. The case where $\tau$ is an $\langmaybe$-function symbol is given. Now suppose the claim holds for $\tau_{1}, \ldots, \tau_{k}$ and consider $\tau(\overline{x}, \overline{y}) = f(\tau_{1}(\overline{x}, \overline{y}), \ldots, \tau_{k}(\overline{x}, \overline{y}))$. Then:
		\begin{align*}
			\tau^{\model{N}}(\overline{d}, \overline{e}) &= f^{\model{N}}(\tau^{\model{N}}_{1}(\overline{d}, \overline{e}), \ldots, \tau^{\model{N}}_{k}(\overline{d}, \overline{e})) \\
			& = f^{\model{N}}(\tau_{1}^{\model{N}}(\overline{a}, \overline{e}), \ldots, \tau^{\model{N}}_{k}(\overline{a}, \overline{e})) = \tau^{\model{N}}(\overline{a}, \overline{e})\\
			&= f^{\model{N}}(\tau_{1}^{\model{M}}(\overline{a}, \overline{e}), \ldots, \tau^{\model{M}}_{k}(\overline{a}, \overline{e}))\\
			& = f^{\model{M}}(\tau_{1}^{\model{M}}(\overline{a}, \overline{e}), \ldots, \tau^{\model{M}}_{k}(\overline{a}, \overline{e})) = \tau^{\model{M}}(\overline{a}, \overline{e})
		\end{align*}
	This proves the claim. Hence, for all atomic $\phi \in \langno_{m+n}$, all $\overline{d} \in D^{m}$ and all $\overline{e} \in M^{n}$:
		\begin{align*}
			\model{N} \models \phi(\overline{d}, \overline{e}) &\text{ iff }\model{N} \models \phi(\overline{a}, \overline{e})\text{ iff }\model{M} \models \phi(\overline{a}, \overline{e})
		\end{align*}
	So for all $d \in D$ we have $a \foi d$ in $\model{N}$ by Lemma \ref{lem:AlternateCharacterisationOfFOI}; and moreover $\model{N} \succno \model{M}$.
\end{proof}
\end{lem}\noindent
Thus armed, I can connect $\foi$ with symmetry in elementary extensions:
\begin{lem}\label{thm:ElementExtensionMain} For any $\langmaybe$-structure $\model{M}$: if $a \symbare b$ in every $\model{N} \succno \model{M}$, then $a \foi b$ in $\model{M}$.
\begin{proof}
	Suppose  $a \symbare b$ in every $\model{N} \succno \model{M}$. Let $D$ be such that $M \cap D = \emptyset$ and $|D| > |\quotient{b}_\model{M}|$. Construct $\model{N}$ as in Lemma \ref{lem:Inflate}, so that $a \foi d$ in $\model{N}$ for all $d \in D$. Since $\model{N} \succno \model{M}$, by assumption there is a symmetry $\pi$ on $\model{N}$ such that $\pi(a) = b$.  So, for every $\phi \in \langno_{2}$, and all $d \in D$, by Lemma \ref{lem:IsomorphismPreservation}:
		\begin{align*}
			\model{N} \models \phi(b, b)\text{ iff }
			\model{N} \models \phi(a, a) \text{ iff }\model{N} \models \phi(a, d) \text{ iff }\model{N} \models \phi(b, \pi(d)) 
		\end{align*}
	Hence $\pi(d) \in \quotient{b}_\model{N}$ for every $d \in D$, by Lemma \ref{lem:AlternateCharacterisationOfFOI}. Since $\pi$ is a bijection, $|D| = |\{\pi(d) \mid d \in D\}| \leq |\quotient{b}_\model{N}|$. If $a \nfoi b$ in $\model{N}$, then $|\quotient{b}_\model{N}| = |\quotient{b}_\model{M}|$, contradicting our choice of $D$; so $a \foi b$ in $\model{N}$. Since $\model{M} \precno \model{N}$ and $\foi$ is universally capturable$\nostroke$ by Lemma \ref{lem:CaptureIndiscernible}, $a \foi b$ in $\model{M}$.
	\end{proof}
\end{lem}\noindent
Theorem \ref{thm:BigMap} (left-diagram) entails us that there is no converse to \ref{thm:ElementExtensionMain} in the general case. However, we do obtain a converse in restricted circumstances:
\begin{lem}When $\langmaybe$ is relational, for any $\langmaybe$-structure $\model{M}$: $a \foi b$ in $\model{M}$ iff $a \symbare b$ in every $\model{N} \succno \model{M}$.
	\begin{proof}
		Immediate from Theorem \ref{thm:BigMap} and Lemmas \ref{lem:CaptureIndiscernible} and \ref{thm:ElementExtensionMain}.
	\end{proof}
\end{lem}\noindent Moreover, we can strengthen Lemma \ref{thm:ElementExtensionMain} in the case of $\symtotal$. 
\begin{lem}\label{thm:ElementExtensionTotal} 
	Let $\model{M}$ be an $\langmaybe$-structure with $a \in M$ and $e \notin M$. Let $\model{N} \succno \model{M}$ be constructed as in Lemma  \ref{lem:Inflate}, so that $N = M \cup \{e\}$ and $a \foi e$ in $\model{N}$. If $a \symtotal b$ in $\model{N}$, then $a \foi b$ in $\model{M}$.
\begin{proof}
 	Suppose $a \symtotal b$ in $\model{N}$, i.e.\ $\pi(a) = b$, $\pi(b) = a$, and $\pi(x) = x$ for all $x \notin \{a, b\}$ is a symmetry on $\model{N}$. In particular, $\pi(e) = e$. Hence, invoking Lemma \ref{lem:IsomorphismPreservation},  for all $\phi \in \langno_{2}$: $\model{M}\models \phi(a, a)$ iff $ \model{N} \models \phi(a, a)$ iff $\model{N} \models \phi(e, a)$ iff $\model{N} \models \phi(e, b)$ iff $\model{N} \models \phi(a, b)$ iff $\model{M} \models \phi(a, b)$. Hence $a \foi b$ in $\model{M}$ by Lemma \ref{lem:AlternateCharacterisationOfFOI}.
\end{proof}
\end{lem}\noindent
However, this strengthening of Lemma \ref{thm:ElementExtensionMain} is limited to the case of $\symtotal$. To see this, let $\model{\TripleRotate}$ comprise two disjoint copies of the complete countable graph, i.e.:
			\begin{align*}
				\TripleRotate &= \mathbb{N}\\
				R^{\model{\TripleRotate}} &= \left\{ \langle m, n\rangle  \in \mathbb{N}^2 \mid m\neq n \text{ and }m + n \text{ is even}\right\}
			\end{align*}
	Whilst $1 \nfoi 2$ in $\model{\TripleRotate}$, we can use Lemma \ref{lem:Inflate} to add a single new element, $e$, such that $1 \foi e$, without disrupting the fact that $1 \sympair 2$. 	Moreover, nothing like Lemma \ref{thm:ElementExtensionMain} holds for relativities: Lemma \ref{lem:RelativityCapturable} tells us that $a \relativetotal b$ in $\model{M}$ iff $a \relativetotal b$ in all $\model{N} \succno \model{M}$. 
	
The results of \S\ref{s:DefinabilityTheory} exhaustively detailed the connections between grades of discernibility and the existence of a symmetry/relativity in \emph{some} elementary extension. The results of this section now exhaustively detail the connections between grades of discernibility and the existence of a symmetry/relativity in \emph{all} elementary extensions. We thus have complete answers to several natural questions concerning the connection between grades of discrimination and elementary extensions.

\section{Concluding remarks}
Several recent technical-cum-philosophical papers have explored some of the grades of discrimination. This paper has pressed forward that technical investigation in many ways. To close, I shall emphasise two.

First, I have introduced \emph{grades of relativity} to the philosophical literature---along with the notion of a near-correspondence, a relativeness correspondence, and a partial relativeness correspondence---and shown that these are the natural $\langno$-analogues of the grades of symmetry. 

Second, I have offered complete answers to the natural questions that arise concerning all twelve grades of discrimination. Indeed, the technical investigation of the grades of discrimination now seems to be complete.\footnote{Huge thanks to \O{}ystein Linnebo and Sean Walsh, whose questions provided much of the original motivation for this paper, and whose subsequent comments were very helpful. Further thanks to Denis Bonnay, Adam Caulton, Fredrik Engstr\"{o}m, Jeffrey Ketland, Richard Pettigrew, a referee for \emph{RSL}, and a referee for this journal.}

\bibliographystyle{jflnat} 
\bibliography{GTA}

\end{document}